\documentclass[12pt,a4paper]{article}
\usepackage{amsmath,amsfonts,amssymb,amsthm,amscd}
\usepackage{hyperref}
\usepackage{braket}
\usepackage{epstopdf}

\usepackage{braket}
\usepackage{xypic}
\usepackage{amsmath}

\usepackage[american]{babel}
\usepackage{color} 
\usepackage{epsfig}
\usepackage{caption}
\usepackage{rotating}
\usepackage{setspace}
\usepackage{fancyhdr}
\usepackage{booktabs}
\usepackage{mathrsfs}
\usepackage{amsthm}
\usepackage{amsmath,amssymb}
\usepackage{enumerate}
\usepackage{tikz}
\usepackage{youngtab}
\usepackage{hyperref}
\usepackage{tensor}
\usepackage{mathabx}
\usepackage{booktabs}
\usepackage{mathrsfs}
\usepackage{subfigure}
\usepackage{amsmath,amssymb}
\usepackage{quiver}

\def\id{{\rm id}}

\renewcommand{\[}{\begin{equation}}
\renewcommand{\]}{\end{equation}}

\newtheorem{thm}{Theorem}[section]
\newtheorem{cor}[thm]{Corollary}
\newtheorem{lem}[thm]{Lemma}
\newtheorem{prop}[thm]{Proposition}

\theoremstyle{definition}
\newtheorem{defn}{Definition}[section]

\theoremstyle{remark}
\newtheorem{oss}{Remark}[section]

\DeclareMathOperator{\Aut}{Aut}

\DeclareMathOperator{\supp}{supp}
\DeclareMathOperator{\SL}{SL}

\DeclareMathOperator{\cohull}{cohull}

\title{Examples of non-amenable, boundary-amenable dynamical systems}
\author{J. Bassi\\ \\Department of Mathematics, University of Tor Vergata, \\Via della Ricerca Scientifica 1, 00133, Roma, Italy\\ \\email: jacopo.bassi@mat.uniroma2.it}

\begin{document}

\maketitle

\begin{abstract}
\noindent Let $\Gamma$ be a discrete countable group with the (AP)-property. It is shown that if $\Gamma$ acts on a countable set $\mathfrak{X}$ in such a way that the infinite intersection of stabilizer subgroups is always trivial, then the induced action of $\Gamma$ on $\partial_\beta \mathfrak{X}$ is topologically amenable. The range of applications include the action of $\Gamma$ on $\partial_\beta (\Gamma / \Lambda)$ for: (i) $\Gamma$ countable hyperbolic torsion-free and $\Lambda$ quasi-isometrically embedded with infinite index, (ii) $\Gamma= \Lambda * \Lambda '$ with $\Lambda$ non-amenable countable, $\Lambda'$ infinite countable and $\Gamma$ with the (AP)-property; moreover this includes the case of actions of groups of automorphisms of a $k$-regular tree with $k \geq 3$ generated by a finite number of Haar-random elements on the Stone-{\v C}ech boundary of the tree.

The techniques involved rely on a study of dynamical properties for actions on non-standard boundaries studied by the author and F. R{\u a}dulescu in previous works.





\end{abstract}
\tableofcontents

\section{Introduction}

Dynamical properties of group actions, and in particular of boundary actions, represent a central topic in modern Analytic Group Theory and certain group properties can be deduced by specializing to particular dynamical systems. For example, exactness of a discrete group $\Gamma$ can be characterized by amenability of the action of $\Gamma$ on either its Stone-{\v C}ech remainder $\partial_\beta \Gamma$ or equivalently its Furstenberg boundary $\partial_F \Gamma$ (\cite{Oz, KaKe}). This latter boundary also contains information about the ideal structure of the reduced group $C^*$-algebra $C^*_\lambda \Gamma$, its simplicity being equivalent to freeness of the action of $\Gamma$ on $\partial_F \Gamma$ (\cite{KaKe}). Equivariant compactifications associated to more general dynamical systems are also of great importance, for example, if $\Gamma$ is biexact (i.e. the left-right action of $\Gamma \times \Gamma$ on $\partial_\beta \Gamma$ is topologically amenable), then the group satisfies the (AO)-property (i.e. the left-right unitary representation of $\Gamma \times \Gamma$ on $\mathbb{B}(l^2 \Gamma)$ becomes tempered after quotienting by the compact operators) and its von Neumann algebra is solid; this was used by N. Ozawa in order to prove that the $II_1$-factor associated to an icc hyperbolic group is prime (\cite{Oz3, AkOs}). Interestingly, the reduced $C^*$-algebra of a biexact group with property T is not nuclear in $K$-theory, hence it does not satisfy the UCT (\cite{Sk}). A much weaker version of biexactness was introduced in \cite{BoIoPe}, where strong rigidity results have been obtained in this much broader context (see also \cite{AmBa} for related results). Boundaries of groups (in the sense of Furstenberg) are related to a conjecture by N. Ozawa about "tightness" of  the embedding of exact $C^*$-algebras in their injective envelopes (\cite{Oz2, KaKe, BaRa3}).

In \cite{BaRa2}, using the deep results contained in \cite{BeKa}, it is observed as boundary amenability of actions on a coset space can be used to deduce information about the ideal structure of the $C^*$-algebra associated to the corresponding quasi-regular representation $\lambda_{\Gamma / \Lambda}$, answering a question posed by Bekka and Kalantar in their work. In fact, if $\Gamma$ is $C^*$-simple, $\Lambda < \Gamma$ has the spectral gap property and the action of $\Gamma$ on $\partial_\beta (\Gamma / \Lambda)$ is topologically amenable, the $C^*$-algebra $C^* (\lambda_{\Gamma / \Lambda} (C^* \Gamma)) \subset \mathbb{B} (l^2 (\Gamma / \Lambda))$ has a unique ideal. It should be emphasized that  the ideal structure of $C^*$-algebras has received a growing interest in the last years. One of the motivations relies on the fact that the structure of simple $C^*$-algebras is very well understood nowadays. To this respect, it is worth mentioning the fundamental results concerning the classification of simple $C ^*$-algebras by means of their Elliott invariant, which is completely settled in the unital case, and for which very important progresses have been made in the non-unital setting (see for example the survey \cite{GoLiNi}). These powerful results can be used to translate dynamical properties in structural properties of $C^*$-algebras in some cases (see for instance \cite{GiPuSk}); in the non-unital setting classification theorems can be used to translate the isomorphism problem for certain subgroups of $\SL(2,\mathbb{R})$ into an isomorphism problem for certain crossed product $C^*$-algebras (\cite{Ba2, Ba1, GoLi}). Leaving the realm of simple $C^*$-algebras, a current line of research concerns the problem of understanding the ideal structure of certain $C^*$-algebras (see for example \cite{ChNe, BrBrLiSc}) and the study of structural properties of $C^*$-algebras with a prescribed set of ideals (\cite{EiReRuSo}). 

In the following we exhibit examples of non-amenable, boundary amenable dynamical systems, some of which provide new examples of $C^*$-algebras associated to quasi-regular representations possessing a unique ideal. The approach used relies on the investigation of dynamical properties of the action of $\Gamma$ on certain spaces, which we will call non-standard boundaries. In the rest of this introduction we recall the original motivation for introducing these dynamical systems and fix some notation.


Let $\Gamma$ be a discrete countable group acting on a countable set $\mathfrak{X}$ by means of bijections, let $\sigma : \Gamma \rightarrow Bij (\mathfrak{X})$ be the given group homomorphism and let $\pi_\sigma : \Gamma \rightarrow U(l^2\mathfrak{X})$ be the associated unitary group representation. Denote by $\mathcal{Q}(\mathfrak{X})$ the \textit{Calkin algebra of $\mathfrak{X}$}, i.e. the $C^*$-algebra $\mathbb{B}(l^2 \mathfrak{X}) /\mathbb{K}(l^2 \mathfrak{X})$ and by $Q : \mathbb{B}(l^2 \mathfrak{X}) \rightarrow \mathcal{Q}(\mathfrak{X})$ the associated quotient $*$-homomorphism. An interesting problem is to understand the relationship between regularity properties of the associated Calkin representation $Q \circ \pi_\sigma : \Gamma \rightarrow \mathcal{Q}(\mathfrak{X})$ and dynamical properties of the associated action $\sigma_\beta$ of $\Gamma$ on $\partial_\beta \mathfrak{X}$, the boundary of the Stone-{\v C}ech compactification of $\mathfrak{X}$. For example, an open problem in the theory is to determine whether the fact that $\Gamma$ satisfies the (AO)-property is enough to guarantee that $\Gamma$ is biexact (\cite{AnDe}). 
Note that the natural generalization of this problem to general actions has a negative answer, since non-exact groups don't act amenably by left translations on their Stone-{\v C}ech boundaries, but the Calkin representations are certainly tempered. 
Motivated by the will to investigate such connection, in \cite{Ba}, based on \cite{BaRa}, it is introduced the notion of non-standard boundaries of a countable set, which are certain extensions of the Stone-{\v C}ech boundary of the given set. The reason for introducing them is that whenever a discrete countable group $\Gamma$ acts on a countable set $\mathfrak{X}$, it is possible to detect particular $\Gamma$-quasi-invariant probability measures on these extensions for which it is possible to compute most vector states for the Koopman representations explicitly and these realize the states on the full group $C^*$-algebra $C^* \Gamma$ coming from the Calkin representation associated to the action of $\Gamma$ on $\mathfrak{X}$. This fact gave in \cite{Ba} the possibility to translate both regularity properties of the Calkin representation $Q \circ \pi_\sigma$ and dynamical properties of the group action $\sigma_\beta$ on $\partial_\beta \mathfrak{X}$ within the common framework of measurable dynamics on non-standard boundaries. The rest of this section contains a recap of some notions and results we will need, plus some results which will be used.
\subsection{Koopman representations}
In this subsection we collect some results concerning quasi-invariant measures and Koopman representations which will be used in the next sections. 

Let $\Gamma$ be a discrete countable group acting on a locally compact Hausdorff space $X$. Recall that a Radon measure $\mu$ on $X$ is called \textit{$\Gamma$-quasi-invariant} if for every $\gamma \in \Gamma$ the pushforward measure $\gamma \mu$ is equivalent to $\mu$ (in the sense of absolute continuity). Now, given such a $\Gamma$-quasi-invariant Radon measure $\mu$ on $X$, there is associated a unitary group representation $\pi_\mu : \Gamma \rightarrow U(L^2 (\mu))$, given by
\begin{equation*}
(\pi_\mu (\gamma) \xi)(x) = \sqrt{\frac{d \gamma \mu}{d \mu}} (x) \xi (\gamma^{-1} x) \qquad \mbox{ for $\mu$-almost every } x \in X,
\end{equation*}
where $d \gamma \mu /d \mu$ is the Radon-Nikodym derivative of $\gamma \mu$ with respect to $\mu$. Combining this group representation with the $GNS$-representation of $C_0 (X)$ gives a covariant pair, hence a $*$-homomorphism $\pi_\mu^\rtimes: C_0 (X) \rtimes \Gamma \rightarrow L^2 (\mu)$, where $C_0 (X) \rtimes \Gamma$ is the full crossed product $C^*$-algebra. From now on we will restrict our attention to the case of actions on compact Hausdorff spaces (rather than locally compact), but the results probably hold more generally.

Given an action of $\Gamma$ on a compact Hausdorff space $X$, there is an interesting connection between the associated topological dynamics and regularity properties of the Koopman representations associated to $\Gamma$-quasi-invariant probability measures on $X$. In order to be more precise, we recall the following
\begin{defn}
The action of $\Gamma$ on $X$ is \textit{topologically amenable} if there is a net of continuous maps $\mu_i : X \rightarrow \mathcal{P}(\Gamma)$ such that
\begin{equation*}
\lim_i \sup_{x \in X} \| \mu_i (\gamma x) - (\gamma \mu_i )(x) \|_1 =0 \qquad \forall \gamma \in \Gamma.
\end{equation*}
\end{defn}
Now, if the action of $\Gamma$ on $X$ is topologically amenable, then for every $\Gamma$-quasi-invariant probability measure $\mu$ on $X$, the associated Koopman representation is \textit{tempered}, i.e. is weakly contained in the left regular representation of $\Gamma$. In order to recover topological amenability from the Koopman representations associated to $\Gamma$-quasi-invariant measures, one needs something stronger than temperedness, namely \textit{Zimmer-amenability} of all these measures (\cite{AnRe} Theorem 3.3.7), i.e. the existence of $\Gamma$-equivariant projections $L^\infty (\mu) \overline{\otimes} l^\infty \Gamma \rightarrow L^\infty (\mu)$, where the action on $L^\infty (\mu) \overline{\otimes} l^\infty \Gamma$ is the diagonal one (we consider the left translation action on $l^\infty \Gamma$); alternatively, for exact groups, Zimmer-amenability has the following characterization:
\begin{lem}[\cite{BuEcWi1} Proposition 3.16]
\label{lem1}
Let $\Gamma$ be a discrete countable exact group acting on a compact Hausdorff space $X$ and let $\mu$ be a $\Gamma$-quasi-invariant probability measure on $X$. Then $\mu$ is Zimmer-amenable if and only if there is a $\Gamma$-equivariant ucp map $l^\infty \Gamma \rightarrow L^\infty (\mu)$.
\end{lem}

Since the action of $\Gamma$ on $X$ is topologically amenable if and only if every $\Gamma$-quasi-invariant probability measure on $X$ is Zimmer-amenable and for every Zimmer-amenable $\Gamma$-quasi-invariant probability measure on $X$ the associated Koopman representation is tempered, it is natural to study the difference between Zimmer-amenability of a measure and temperedness of the associated Koopman representation; it turns out that the missing ingredient is pair-amenability of the measure, as proved in \cite{AnDe2} for the case of second countable locally compact groups acting on standard spaces. We collect this in the following
\begin{thm}[\cite{Ba} Corollary 1.4]
The action of $\Gamma$ on the compact Hausdorff space $X$ is topologically amenable if and only if for every $\Gamma$-quasi-invariant probability measure $\mu$ on $X$ we have
\begin{itemize}
\item[(i)] the associated Koopman representation $\pi_\mu$ is tempered,
\item[(ii)] $\mu$ is pair-amenable, i.e. there is an equivariant norm-one projection $L^\infty (\mu) \overline{\otimes} L^\infty (\mu) \rightarrow L^\infty (\mu)$.
\end{itemize}
\end{thm}
We end this section by collecting two results concerning Koopman representations which will be used in the next sections, where we want to construct inverse limits associated to Koopman representations.

\begin{prop}
\label{prop1}
Let $\Gamma$ be a discrete group acting on a locally compact Hausdorff space $X$ and let $\mu$, $\nu$ be $\Gamma$-quasi-invariant probability measures on $X$ with $\mu \precsim \nu$. There is an isometry $U : L^2 (\mu) \rightarrow L^2 (\nu)$ such that $U(\pi_\mu^{\rtimes} (f\gamma) \xi) = \pi_\nu^{\rtimes} (f \gamma) (U \xi)$ for every $\xi \in L^2 (\mu)$, $f \in C_0(X)$ and $\gamma \in \Gamma$.
\end{prop}
\proof Let $U: L^2 (\mu) \rightarrow L^2 (\nu)$ be given by $U\xi (x) = \sqrt{d\mu /d \nu} (x) \xi (x)$. This is an isometry and for $f \in C_0(X)$, $\gamma \in \Gamma$ and $\xi \in L^2 (\mu)$ we have
\begin{equation*}
(U \pi_\mu^\rtimes (f \gamma) \xi)(x)=\sqrt{\frac{d\mu}{d\nu}} (x) f(x) \sqrt{\frac{d \gamma \mu}{d\mu}} (x) \xi (\gamma^{-1} x)=f(x) \sqrt{\frac{d\mu}{d\nu}} (x)  \sqrt{\frac{d \gamma \mu}{d\mu}} (x) \xi (\gamma^{-1} x),
\end{equation*}
\begin{equation*}
(\pi_\nu^\rtimes (f \gamma) (U\xi) )(x)= f(x) \sqrt{\frac{d\gamma \nu}{d \nu}} (x) \sqrt{\frac{d\mu}{d\nu}}(\gamma^{-1} x) \xi (\gamma^{-1} x)
\end{equation*}
For $\nu$-almost every $x \in X$. Hence we need to check that 
\begin{equation*}
\sqrt{\frac{d\mu}{d\nu}} (x)  \sqrt{\frac{d \gamma \mu}{d\mu}} (x) =  \sqrt{\frac{d\gamma \nu}{d \nu}} (x) \sqrt{\frac{d\mu}{d\nu}}(\gamma^{-1} x)
\end{equation*}
 for $\nu$-almost every $x \in X$. Now (see \cite{Fur} 3.1) for every $g \in C_0 (X)$ we have
\begin{equation*}
\begin{split}\int g(x) \frac{d \gamma \mu}{d \gamma \nu} (x) d\gamma \nu (x) &= \int g (x) d \gamma \mu (x) = \int g (\gamma x) d\mu (x) = \int g(\gamma x) \frac{d\mu}{d\nu} (x) d\nu (x) \\&= \int g(x) \frac{d\mu}{d\nu} (\gamma^{-1} x) d\gamma \nu (x), \end{split}
\end{equation*}
from which it follows that $\frac{d \gamma \mu}{d \gamma \nu}= \frac{d\mu}{d\nu} \circ \gamma^{-1}$ for $\gamma \nu$-almost every $x$ (and so for $\nu$-almost every $x$). The result follows since $\frac{d\mu}{d\nu}   \frac{d \gamma \mu}{d\mu} =  \frac{d \gamma \mu}{d\nu} = \frac{ d\gamma \nu}{d \nu} \frac{d\gamma \mu}{d\gamma \nu}$ $\nu$-almost everywhere. $\Box$\\

The above result clearly has the following counterpart at the level of covariant representations

\begin{cor}
\label{cor1}
Let $\Gamma$ be a discrete group acting on a locally compact Hausdorff space $X$ and $\mu$, $\nu$ be $\Gamma$-quasi-invariant probability measures on $X$ with $\mu \precsim \nu$. There is a surjective $\Gamma$-equivariant $*$-homomorphism $\pi_{\mu, \nu}^\rtimes : \pi_\nu^\rtimes (C_0 (X) \rtimes \Gamma) \rightarrow \pi_\mu^\rtimes (C_0 (X) \rtimes \Gamma)$ such that $\pi_\mu^\rtimes = \pi_{\mu,\nu}^\rtimes \circ \pi_\nu^\rtimes$.
\end{cor}
\proof Proposition \ref{prop1} guarantees that $\ker (\pi_\nu^\rtimes) \subset \ker (\pi_\mu^\rtimes)$. The result follows. $\Box$\\

\subsection{Dynamics on non-standard boundaries}
In this subsection we mainly fix some notation and recall some results from \cite{Ba} which will be used in the following. In the appendix the interested reader can find sketches of proofs of these results.

Let $\Gamma$ be a discrete countable group acting on a countable set $\mathfrak{X}$. As stated above, the necessity to introduce certain extensions of the Stone-{\v C}ech boundary of $\mathfrak{X}$ is motivated by the fact that certain regularity properties of the Calkin representation are encoded in regularity properties of certain $\Gamma$-quasi-invariant probability measures on these extensions. Based on \cite{Ba}, we now recall the relevant notions and results.

Let $\omega \in \partial_\beta \mathbb{N}$ be a free ultrafilter on the natural numbers and consider the $C^*$-algebraic ultraproduct $l^\infty \mathfrak{X}_\omega$, which is naturally endowed with an action of $\Gamma$ (componentwise). We recall that it is defined as the quotient $\prod_\mathbb{N} l^\infty \mathfrak{X} /I_\omega$, where $I_\omega=\{ (f_n) \in \prod_\mathbb{N} \; | \; \lim_{n \rightarrow \omega} \| f_n\| =0\}$. Given $(f_n) \in \prod_\mathbb{N} l^\infty \mathfrak{X}$ we will denote by $(f_n)^\bullet$ the corresponding class in $l^\infty \mathfrak{X}_\omega$, $(f_n)^\bullet = (f_n) + I_\omega$. Since there is a $\Gamma$-equivariant embedding $l^\infty \mathfrak{X} \rightarrow l^\infty \mathfrak{X}_\omega$, there is a $\Gamma$-equivariant continuous surjection $\psi_\omega : \sigma (l^\infty \mathfrak{X}_\omega) \rightarrow \Delta_\beta \mathfrak{X}$, we call such surjection a \textit{non-standard-bounday map}. The \textit{non-standard boundary} associated to $\omega$ is defined to be $\partial_{\beta,\omega} \mathfrak{X} := \psi_\omega^{-1} (\partial_\beta \mathfrak{X})$. We recall from \cite{Ba} Lemma 2.1 that $C(\partial_{\beta, \omega} \mathfrak{X})$ is $\Gamma$-equivariantly $*$-isomorphic to $l^\infty \mathfrak{X}_\omega / \mathcal{C}_{0, \omega} (\mathfrak{X})$, where the ideal $\mathcal{C}_{0,\omega} (\mathfrak{X})$ is obtained as the closure of the union $\bigcup_k \mathcal{C}_{k, \omega}$, where $\mathcal{C}_{k,\omega} (\mathfrak{X}) =\{ (f_n)^\bullet \in l^\infty \mathfrak{X}_\omega \; | \; \supp (f_n) \subset \mathfrak{X}_k \; \forall n \in \mathbb{N}\}$ and the sets $\mathfrak{X}_k \subset \mathfrak{X}$ are a countable increasing family of finite subsets of $\mathfrak{X}$ whose union is $\mathfrak{X}$. \\
Let $(\tilde{\xi}_n)$ be a sequence of norm-one vectors in $l^2 \mathfrak{X}$ which go to zero weakly with respect to $\omega$ and let $(\alpha_i)$ be a sequence of strictly positive real numbers such that $\sum_{i,j} \alpha_i \alpha_j < \infty$. Consider an enumeration of $\Gamma=(\gamma_i)$, with $\gamma_0 = \id$ and define new sequences $(\hat{\xi}_n)$ by $\hat{\xi}_n = \sum_i \alpha_i |\tilde{\xi}_n| \circ \gamma_i^{-1}$ and $(\xi_n)$ by $\xi_n = \hat{\xi}_n /\| \hat{\xi}_n\|_2$. In this way we obtain a $\Gamma$-quasi-invariant probability measure on $\partial_{\beta, \omega} \mathfrak{X}$ through the formula $f \mapsto \lim_{n \rightarrow \omega} \braket{F_n \xi_n, \xi_n}_{l^2 \mathfrak{X}}$, where $(F_n)^\bullet$ is any representative of $f$ in $l^\infty \mathfrak{X}_\omega$, $f= (F_n)^\bullet + \mathcal{C}_{0,\omega}(\mathfrak{X})$. In virtue of the next result, we refer to measures obtained in this way as to \textit{Calkin measures} and denote this set by $\mathcal{P}_{\mathcal{Q}} (\partial_{\beta, \omega} \mathfrak{X})$ (note that this set was denoted by $Q-\tilde{\mathcal{P}}^\infty (\mathfrak{X})_\omega$ in \cite{Ba}). Their role in this setting is encoded in the following
\begin{prop}
\label{prop1.5}
Let $\Gamma$ be a discrete countable group and $\sigma$ an action of $\Gamma$ on a countable set $\mathfrak{X}$, $u_\sigma$ the associated unitary representation on $l^2 \mathfrak{X}$. Let $C^*_\tau \Gamma$ be a $C^*$-completion of the group algebra $\mathbb{C} [\Gamma]$. The Calkin representation $Q \circ u_\sigma : C^* \Gamma \rightarrow  \mathcal{Q}(\mathfrak{X})$ factors through the canonical surjection $\tau: C^* \Gamma \rightarrow C^*_\tau \Gamma$ if and only if for every $\omega \in \partial_\beta \mathbb{N}$ and any $\mu \in \mathcal{P}_\mathcal{Q} (\partial_{\beta, \omega} \mathfrak{X})$ the associated Koopman representation $\pi_\mu : C^* \Gamma \rightarrow \mathbb{B}(L^2 (\mu))$ factors through $\tau$.
\end{prop}
A particular subset of Calkin measures, which will be relevant in the following, are the ones obtained by choosing the sequence $(\tilde{\xi}_n)$ as $\tilde{\xi}_n=\delta_{x_n}$ for a certain sequence $(x_n)$ in $\mathfrak{X}$. As in \cite{Ba} we refer to this set as the set of \textit{uniform approximants of non-standard points} and denote it by $\widetilde{NS}_\omega (\mathfrak{X})$. The terminology is reminiscent of the fact that the set of states on $l^\infty \mathfrak{X}_\omega$ of the form $f \mapsto \lim_{n\rightarrow \omega} f_n (x_n)$ with its $w^*$-topology is homeomorphic to the non-standard version of $\mathfrak{X}$, $^\omega \mathfrak{X}$ with its discrete topology. As observed in \cite{Ba} Proposition 2.11 every $\mu \in \cohull{\widetilde{NS}_\omega (\mathfrak{X})}$ is automatically pair-amenable.

An interesting fact about the topological dynamics on non-standard boundaries is that the non-standard boundary maps are always amenable, in the sense that if the action of $\Gamma$ on $\partial_{\beta, \omega} \mathfrak{X}$ is topologically amenable for some $\omega \in \partial_\beta \mathfrak{X}$, then the action of $\Gamma$ on $\partial_\beta \mathfrak{X}$ is topologically amenable as well. This is indeed a consequence of the following Proposition, which requires some notation (see \cite{Ba} Section 2 for a more detailed treatment). Let $\omega \in \partial_\beta \mathbb{N}$, associated to a given enumeration of $\mathfrak{X}=(\bar{x}_i)$ we can consider for every $k \in \mathbb{N}$ the projection $P_k \in \mathbb{B}(l^2 \mathfrak{X})$ on the linear span of $\{\delta_{\bar{x}_0},...,\delta_{\bar{x}_k}\}$ and, after embedding $P_k$ in $\mathbb{B}(l^2 \mathfrak{X})_\omega$, the $C^*$-algebra $M_{k+1,\omega}:= P_k \mathbb{B}(l^2 \mathfrak{X})_\omega P_k$. Let then $\mathcal{K}_\omega(l^2 \mathfrak{X}) := \overline{\bigcup_k M_{k,\omega}}$ and let $\mathcal{B}_\omega (l^2 \mathfrak{X})$ be the maximal unital $C^*$-subalgebra of $\mathbb{B}(l^2 \mathfrak{X})_\omega$ containing $\mathcal{K}_\omega(l^2 \mathfrak{X})$ as an ideal.
\begin{prop}[\cite{Ba} Proposition 2.3]
\label{prop1.6}
Let $\mathfrak{X}$ be a countable set endowed with an action $\sigma$ of a countable exact group $\Gamma$ and denote by $u_\sigma$ the associated unitary representation on $l^2 \mathfrak{X}$. Let $\omega \in \partial_\beta \mathbb{N}$ and denote by $u_{\sigma, \omega}$ the unitary group representation of $\Gamma$ in $\mathcal{B}_\omega (l^2 \mathfrak{X})$ given by $u_{\sigma, \omega} (\gamma)= (u_\sigma (\gamma))^\bullet$. The following are equivalent:
\begin{itemize}
\item[(i)] The action of $\Gamma$ on $\partial_{\beta, \omega} \mathfrak{X}$ is topologically amenable,
\item[(ii)] there is a nuclear ucp map $\phi: C^*_\lambda \Gamma \rightarrow \mathcal{B}_\omega (l^2 \mathfrak{X})$ such that $\phi (\lambda (\gamma)) - u_{\sigma, \omega} (\gamma) \in \mathcal{K}_\omega (l^2 \mathfrak{X})$ for every $\gamma \in \Gamma$,
\item[(iii)] the action of $\Gamma$ on $\partial_\beta \mathfrak{X}$ is topologically amenable.
\end{itemize}
\end{prop}

Note that the above result is really a "non-standard version" of \cite{Oz4} Proposition 4.1.\\
 In order to apply the above Proposition in what follows, we observe in the following Lemma that the covariant representations of $(l^\infty \mathfrak{X}_\omega, \Gamma)$ which arise from Calkin measures always give rise to Koopman representations which factor through the homomorphism associated to the covariant representation of $(l^\infty \mathfrak{X}_\omega, \Gamma)$ in $\mathcal{B}(l^2 \mathfrak{X})_\omega$.
 
 \begin{lem}
\label{lemfac}
Let $\omega \in \partial_\beta \mathbb{N}$ and let $\mu$ be a Calkin measure on $\sigma (l^\infty \mathfrak{X}_\omega)$. The Koopman representation $\pi_\mu^\rtimes : l^\infty \mathfrak{X}_\omega \rtimes \Gamma \rightarrow \mathbb{B}(L^2(\mu))$ factors through the covariant representation $u_{\sigma, \omega}^\rtimes : l^\infty X_\omega \rtimes \Gamma \rightarrow \mathcal{B}_\omega(l^2 \mathfrak{X})$, i.e. the following diagram is commutative:
\[\begin{tikzcd}
	{l^\infty \mathfrak{X}_\omega \rtimes \Gamma} && {\mathbb{B}(L^2 (\mu))} \\
	\\
	& {C^*(l^\infty \mathfrak{X}_\omega, u_{\sigma, \omega}(\Gamma))} \\
	& {\mathcal{B}_\omega(l^2 \mathfrak{X})}
	\arrow["{\pi_\mu^\rtimes}", from=1-1, to=1-3]
	\arrow["u_{\sigma, \omega}^\rtimes"{description}, from=1-1, to=3-2]
	\arrow["\exists"{description}, from=3-2, to=1-3]
	\arrow[hook, from=3-2, to=4-2]
\end{tikzcd}\]
\end{lem}
\proof As observed in the proof of \cite{Ba} Propostion 2.13, the vectors of the form $[(\eta_n/\tilde{\xi}_n)]$ with $\eta_n=\sum_{i=0}^N \alpha_i f_n (\xi_n \circ s_i^{-1})$ for every $n$ and $(f_n)^\bullet \in l^\infty \mathfrak{X}_\omega$ are dense in $L^2(\mu)$. Elements of the form $\sum_{j=1}^k h_j \gamma_j \in l^\infty \mathfrak{X}_\omega \rtimes \Gamma$, with $h_j \in l^\infty \mathfrak{X}_\omega$, $\gamma_j \in \Gamma$, $j=1,...,k$, act on these vectors as $\pi_\mu (\sum_{j=1}^k h_j \gamma_j) [(\eta_n /\tilde{\xi}_n)]= [((\sum_{i=0}^N \sum_{j=1}^k \alpha_i h_{j,n} (\gamma_j f_n) (\xi_n \circ s_i^{-1})/\tilde{\xi}_n)]$. In particular, $\|[((\sum_{i=0}^N \sum_{j=1}^k \alpha_i h_{j,n} (\gamma_j f_n) (\xi_n \circ s_i^{-1})/\tilde{\xi}_n)]\|_{L^2(\mu)} = \lim_{n \rightarrow \omega} \|\sum_{i=0}^N \sum_{j=1}^k \alpha_i h_{j,n} (\gamma_j f_n) (\xi_n \circ s_i^{-1})\|_{l^2 \mathfrak{X}}$, from which it follows that $\|\pi_\mu (a)\| \leq \|u_{\sigma , \omega}^\rtimes (a)\|$ for a dense set of elements in $l^\infty \mathfrak{X}_\omega \rtimes \Gamma$. $\Box$\\

\section{Non-trivial ultrafilters}
In this section we want to specialize our approach to the case of non-trivial ultrafilters on $\mathbb{N}$. Recall that $\omega \in \partial_\beta \mathbb{N}$ is called non-trivial if there is a sequence of sets $(A_n) \subset \omega$ such that $\cap_n A_n = \emptyset$. Note that $\partial_\beta \mathbb{N}$ contains non-trivial ultrafilters, indeed one can take any maximal extension of the filter $\alpha=\{ A \subset \mathbb{N} \; | \; |A^c | <\infty\}$. If $\omega \in \partial_\beta \mathbb{N}$ is non-trivial, then there is a strictly decreasing sequence $B_0 \supset B_1 \supset B_2 \supset ...$ in $\omega$ with $\cap_n B_n = \emptyset$. For example one could take $B_n = \cap_{i=0}^{j_n} A_i$, where $j_0 =1$ and, for $n \geq 1$, $j_n$ is given by $j(n)= \min \{ k > j_{n-1} \; | \; \cap_{i=0}^{j_{n}} A_i \subsetneq \cap_{i=0}^{j_{n-1}} A_i\}$. Given $\omega \in \partial_\beta \mathbb{N}$ and a countable set $\mathfrak{X}$, we denote by ${}^\omega \mathfrak{X}$ the non-standard model of $\mathfrak{X}$ with respect to $\omega$, that is the quotient $\{(x_n) \subset \mathfrak{X}\} /\sim_\omega$, where $(x_n) \sim_\omega (y_n)$ if there is $A \in \omega$ such that $x_n =y_n$ for every $n \in A$, and by ${}^\omega \mathfrak{X}_\infty$ the part at infinity of the non-standard model of $\mathfrak{X}$, that is the set of classes $[x_n]_\omega$ such that $x_n \rightarrow \infty$ with respect to $\omega$ (this definition does not depend on the representative).

Before stating the results of this section we give a brief overview. Let $\Gamma$ be a discrete countable group acting on a countable set $\mathfrak{X}$. First we show that given a non-trivial free ultrafilter $\omega$ on $\mathbb{N}$, the part at infinity of the non-standard model ${}^\omega \mathfrak{X}_\infty$ (with its discrete topology) is $\Gamma$-equivariantly embedded as a discrete, open dense subset of $\partial_{\beta, \omega} \mathfrak{X}$ (Proposition \ref{prop2.1}). This fact will be used later on in order to guarantee that if a function $f \in C(\partial_{\beta, \omega} \mathfrak{X})$ vanishes on ${}^\omega \mathfrak{X}_\infty$, then it vanishes everywhere on $\partial_{\beta, \omega} \mathfrak{X}$. Next we show that, using the same procedure outlined for example in the proof of \cite{BaRa} Theorem 2.7 , given a point $[x_n]_\omega \in {}^\omega \mathfrak{X}_\infty$, we can construct a sequence $\xi=(\xi_n) \subset l^2 \mathfrak{X}$ for which the measure defined by $C(\partial_{\beta, \omega} \mathfrak{X}) \rightarrow \mathbb{C}$, $f \mapsto \lim_{n \rightarrow \omega} \braket{F_n \xi_n, \xi_n}$, where $F=(F_n)^\bullet$ is a representative of $f$ in $l^\infty \mathfrak{X}_\omega$, is a $\Gamma$-quasi-invariant probability measure supported on ${}^\omega \mathfrak{X}_\infty$ (Proposition \ref{prop2.2}). Then we use measures constructed in this way in order to construct an inverse limit $\pi$ of Koopman representations of $C^* (l^\infty \mathfrak{X}_\omega, u_{\sigma, \omega} (\Gamma)) \subset \mathcal{B}_\omega(l^2 \mathfrak{X})$ (in principle these are representations of $l^\infty \Gamma \rtimes \Gamma$, but we can employ Lemma \ref{lemfac}) which allow, in the case the stabilizers of point in ${}^\omega \mathfrak{X}_\infty$ are amenable, the construction of a nuclear ucp map from $C^*_\lambda \Gamma$ to $\pi (C^* (l^\infty \mathfrak{X}_\omega, u_{\sigma, \omega} (\Gamma)))$ which gives , by Choi-Effros lifting (\cite{ChEf} Theorem 3.10), a nuclear ucp map as the one appearing in Proposition \ref{prop5} (ii), except for the fact that we do not have control on the difference $\phi (\lambda (\gamma)) - u_{\sigma, \omega} (\gamma)$ for $\gamma \in \Gamma$. In the final part of this section we see how the hypothesis that the stabilizers of points in ${}^\omega \mathfrak{X}_\infty$ are actually trivial (rather than just amenable), together with the (AP)-property imply that $\ker \pi \subset \mathcal{K}_\omega(l^2 \mathfrak{X})$ and so that the action of $\Gamma$ on $\partial_{\beta} \mathfrak{X}$ is topologically amenable.

\begin{prop}
\label{prop2.1}
Let $\Gamma$ be a discrete countable group acting on a countable set $\mathfrak{X}$. Let $\omega \in \partial_\beta \mathbb{N}$ be non-trivial. Then ${}^\omega \mathfrak{X}_\infty$ is a discrete dense subset of $\partial_{\beta, \omega} \mathfrak{X}$. In particular it is open and $\partial_{\beta, \omega} \mathfrak{X}$ is a $\Gamma$-equivariant compactification of ${}^\omega \mathfrak{X}_\infty$.
\end{prop}
\proof Note that every $[x_n]_\omega \in {}^\omega \mathfrak{X}_\infty$ gives rise to a character on $l^\infty \mathfrak{X}_\omega$ given by $(f_n)^\bullet \mapsto \lim_{n \rightarrow \omega} \braket{f_n  \delta_{x_n}, \delta_{x_n}}$, which factors through $\mathcal{C}_{0,\omega} (\mathfrak{X})$. The embedding ${}^\omega \mathfrak{X}_\infty$ is given by $[x_n]_\omega \mapsto \{f \mapsto \lim_{n \rightarrow \omega} \braket{F_n \delta_{x_n}, \delta_{x_n}}\; | \;  (F_n)^\bullet + \mathcal{C}_{0,\omega}(\mathfrak{X}) = f, \; (F_n)^\bullet \in l^\infty \mathfrak{X}_\omega \}$. This embedding is $\Gamma$-equivariant. In order to check that the image is a discrete subset of $\partial_{\beta, \omega} \mathfrak{X}$, we proceed as in \cite{Ba} Remark 2.1. Hence let $[x_n]_\omega$ be given and view it as a character on $C(\partial_{\beta, \omega} \mathfrak{X})$. Consider the open set around $[x_n]_\omega$ given by $U=\{ y \in \partial_{\beta, \omega} \mathfrak{X} \; | \; | [x_n]_\omega ((\delta_{x_n})^\bullet + \mathcal{C}_{0,\omega} (\mathfrak{X})) - y((\delta_{x_n})^\bullet + \mathcal{C}_{0,\omega} (\mathfrak{X}))| < 1/2\}$. Let $[y_n]_\omega \neq [x_n]_\omega$ and suppose $[y_n]_\omega \in U$ and note that for every $n$ we have $(y_n)(\delta_{x_n}) \in \{0,1\}$. Hence $[y_n]_\omega ((\delta_{x_n})^\bullet + \mathcal{C}_{0,\omega} (\mathfrak{X})) \in \{0,1\}$. Now, if $[y_n]_\omega ((\delta_{x_n})^\bullet + \mathcal{C}_{0,\omega} (\mathfrak{X})) =\lim_{n\rightarrow \omega} \braket{ \delta_{x_n} \delta_{y_n}, \delta_{y_n}}=1$, there is $A \in \omega$ such that $y_n = x_n$ for every $n \in A$. Hence, if $[y_n]_\omega \in U$, then $[y_n]_\omega = [x_n]_\omega$ and ${}^\omega \mathfrak{X}_\infty$ is a discrete subset of $\partial_{\beta, \omega} \mathfrak{X}$. Let now $(A_n)$ be a strictly decreasing sequence of sets in $\omega$ with empty intersection. Let $f=(f_n)^\bullet \in l^\infty \mathfrak{X}_\omega$ be such that the following holds: there is $c>0$ such that for every $k \in \mathbb{N}$ we have $\| f- P_k f P_k\| >c$ (see the paragraph before Proposition \ref{prop1.6} for the definition of the projections $P_k$). For every $k \in \mathbb{N}$ let $\mathfrak{X}_k = \{\bar{x}_0, \bar{x}_1,..., \bar{x}_k\}$. Let $g^{(1)} = f- P_1 f P_1 = (f_n - P_1 f_n P_1)^\bullet$ and $g_n^{(1)} := f_n - P_1 f_n P_1$ for every $n \in \mathbb{N}$. There is a sequence $(x^{(1)}_n) \subset \mathfrak{X} \backslash \mathfrak{X}_1$ such that $\| g^{(1)}\|=\lim_{n \rightarrow \omega} | g^{(1)}_n (x^{(1)}_n)|>c$. Hence there is $B_1 \in \omega$ such that $| g^{(1)}_n (x^{(1)}_n)| >c$ for every $n \in B_1$. Let $A_k$ be a strictly decreasing sequence of subsets of $\omega$ with empty intersection and take $C_1 = A_1 \cap B_1$. Suppose now that, given $2\leq  k \in \mathbb{N}$, there are a sequence $(x^{(k)}_n) \subset \mathfrak{X} \backslash \mathfrak{X}_k$ and a set $C_k \in \omega$, with $C_k \subsetneq C_{k-1}$, satisfying $|g^{(k)}_n (x^{(k)}_n)| >c$ for every $n \in C_k$, where $g^{(k)}_n = f_n - P_k f_n P_k$. By hypothesis, given $g^{(k+1)}=(g^{(k+1)}_n)^\bullet$, where $g^{(k+1)}_n = f_n - P_{k+1} f_n P_{k+1} \in l^\infty \mathfrak{X}$, there is a sequence $(x^{(k+1)}_n) \subset \mathfrak{X}\backslash \mathfrak{X}_{k+1}$ such that $\|g^{(k+1)}\| = \lim_{n \rightarrow \omega} | g^{(k+1)}_n (x_n)| >c$. In particular, there is $B_{k+1} \in \omega$ such that $|g^{(k+1)}_n (x^{(k+1)}_n)| > c$ for every $n \in B_{k+1}$. Take $C_{k+1} = B_{k+1} \cap A_{k+1}$. Define the following sequence $(x_n) \subset \mathfrak{X}$: 
\begin{equation*}
x_n =\begin{cases}  \bar{x}_0 & \mbox{ for } n \in \mathbb{N} \backslash C_1\\
 x^{k}_n & \mbox{ for } n \in C_k \backslash C_{k-1} \mbox{ for some } k.
 \end{cases}
\end{equation*}
Note that since $\cap_k C_k = \emptyset$ the above really defines a sequence. 
The sequence $(x_n)$ goes to $\infty$ with respect to $\omega$ and satisfies $\lim_{n \rightarrow \omega} |f_n (x_n)| >c$. We have just shown that if there is no sequence $(x_n) \subset \mathfrak{X}$ going to infinity with respect to $\omega$ such that $\lim_{n \rightarrow \omega} |f_n (x_n)| >0$, then for every $\epsilon >0$ there is $k \in \mathbb{N}$ such that $\| f- P_k f P_k\| < \epsilon$, which gives $f \in \mathcal{C}_{0,\omega} (\mathfrak{X})$. Hence for every non-zero function $f \in C(\partial_{\beta, \omega} \mathfrak{X})$ there is a character given by such a sequence $(x_n)$ on which $f$ does not vanish, which, since $\partial_{\beta, \omega} \mathfrak{X}$ is Hausdorff, implies that ${}^\omega \mathfrak{X}_\infty$ is dense in $\partial_{\beta, \omega} \mathfrak{X}$. Now, it is a general fact that the points in a discrete dense subset $Y$ of a Hausdorff space $X$ are open, indeed suppose this is not the case, then there is an element $y \in Y$ such that every open neighborhood of $y$ contains a point in $X \backslash Y$; choose such an open set such that it does not contain any other point in $Y$ except $y$. This contains a point $z \in X \backslash Y$ and so an open set around $z$ which does not intersect $Y$, which is impossible since $Y$ is dense. $\Box$\\

Given a uniformly bounded sequence of vectors $\xi= (\xi_n)$ in $l^2 \mathfrak{X}$ we will denote by $\mu_\xi$ the finite Radon measure on $\sigma (l^\infty \mathfrak{X}_\omega)$ given by $(f_n)^\bullet \mapsto \lim_{n \rightarrow \omega} \braket{ f_n \xi_n ,\xi_n}$.

\begin{lem}
\label{lemnew1}
Let $\Gamma$ be a discrete countable group acting on a countable set $\mathfrak{X}$. For every uniformly bounded sequence $\xi=(\xi_n)$ in $l^2 \mathfrak{X}$ such that $\lim_{n \rightarrow \omega} \|\xi_n\|_2 >0$ we have $\mu_{(\xi_n /\| \xi_n\|_2)} = \mu_{\xi/\lim_{n \rightarrow \omega} \|\xi_n\|_2}$.
\end{lem}
\proof Let $f =(f_n)^\bullet \in l^\infty \mathfrak{X}_\omega$ and $\epsilon >0$ be given, choose $\delta >0$ such that $0 < \delta < \lim_{n \rightarrow \omega} \|\xi_n\|_2 /2$ and $2\delta + 3 \delta \| f\|_\infty (\lim_{n \rightarrow \omega} \| \xi_n \|_2)^{-1}   < \epsilon$. Pick $A \in \omega$ such that 
 \begin{equation*}
 \begin{split}
 \max\{& | \mu_{(\xi_n /\| \xi_n\|_2)} (f) - \braket{f_k \xi_k /\|\xi_k\|_2 , \xi_k / \|\xi_k\|_2}|, |\mu_{\xi/\lim_{n \rightarrow \omega} \|\xi_n\|_2} (f) - \braket{f_k \xi_k /\lim_{n \rightarrow \omega} \|\xi_n\|_2, \xi_k /\lim_{n \rightarrow \omega} \|\xi_n\|_2}|, \\ &|\lim_{n\rightarrow \omega} \|\xi_n\|_2 - \|\xi_k\|_2| \} < \delta
 \end{split}
 \end{equation*}
 for every $k \in A$. For such $k$'s we have 
 \begin{equation*}
 \begin{split}
 |\mu_{(\xi_n /\| \xi_n\|_2)} (f) &- \mu_{\xi/\lim_{n \rightarrow \omega} \|\xi_n\|_2} (f) | \leq |\mu_{(\xi_n /\| \xi_n\|_2)} (f) - \braket{f_k \xi_k /\|\xi_k\|_2 , \xi_k / \|\xi_k\|_2}| \\ & + |\mu_{\xi/\lim_{n \rightarrow \omega} \|\xi_n\|_2} (f) - \braket{f_k \xi_k /\lim_{n \rightarrow \omega} \|\xi_n\|_2, \xi_k /\lim_{n \rightarrow \omega} \|\xi_n\|_2}| \\ &+ |\braket{f_k \xi_k /\|\xi_k\|_2 , \xi_k / \|\xi_k\|_2} - \braket{f_k \xi_k /\lim_{n \rightarrow \omega} \|\xi_n\|_2 , \xi_k /\lim_{n \rightarrow \omega} \|\xi_n\|_2}|\\
 & < 2 \delta + |\braket{f_k \xi_k /\|\xi_k\|_2 , \xi_k / \|\xi_k\|_2} - \braket{f_k \xi_k /\lim_{n \rightarrow \omega} \|\xi_n\|_2 , \xi_k /\|\xi_k\|_2}| \\ & + |\braket{f_k \xi_k /\lim_{n \rightarrow \omega} \|\xi_n\|_2 , \xi_k /\|\xi_k\|_2} - \braket{f_k \xi_k /\lim_{n \rightarrow \omega} \|\xi_n\|_2 , \xi_k /\lim_{n \rightarrow \omega} \|\xi_n\|_2}|\\
 & \leq 2 \delta + \|f\|_\infty \|\xi_k\|_2 |\|\xi_k\|_2 - \lim_{n \rightarrow \omega} \|\xi_n\|_2| \|\xi_k\|^{-1}_2 (\lim_{n \rightarrow \omega} \|\xi_n \|_2)^{-1}\\
 & + \|f\|_\infty \|\xi_k\|_2 (\lim_{n \rightarrow \omega} \|\xi_n \|_2)^{-1} \|\xi_k\|_2 | \|\xi_k \|_2 - \lim_{n \rightarrow \omega} \|\xi_n \|_2| \| \xi_k \|^{-1}_2 (\lim_{n \rightarrow \omega} \|\xi_n \|_2)^{-1}\\
 & < 2\delta + \|f\|_\infty \delta ((\lim_{n\rightarrow \omega} \|\xi_n \|_2)^{-1} + \|f\|_\infty \delta \|\xi_k \|_2 (\lim_{ n \rightarrow \omega} \|\xi_n\|_2)^{-2})\\
 & \leq 2 \delta + \delta \left(\| f\|_\infty (\lim_{n \rightarrow \omega} \| \xi_n \|_2)^{-1} (1 + (\lim_{n \rightarrow \omega} \| \xi_n \|_2 + \delta)(\lim_{n \rightarrow \omega} \| \xi_n \|_2)^{-1} \right)\\
 & < \epsilon.
 \end{split}
 \end{equation*}
 The result follows. $\Box$

\begin{prop}
\label{prop2.2}
Let $\Gamma$ be a discrete countable group acting on a countable set $\mathfrak{X}$. Let $(x_n)$ be a sequence in $\mathfrak{X}$ such that $\lim_{n \rightarrow \omega} x_n = \infty$ and let $\Lambda= Stab ([x_n]_\omega) \subset \Gamma$. Then for every sequence $(\alpha_i) \subset l^1 \mathbb{N}$ of strictly positive real numbers, the associated probability measure $\mu_{(\sum_i \alpha_i \delta_{(x_n)} \circ \gamma^{-1}_i)}$ on $\partial_{\beta, \omega} \mathfrak{X}$ coincides with the $\Gamma$-quasi inviariant probability measure $\tilde{\mu}_{(\alpha_i), [x_n]_\omega}$ on $\partial_{\beta, \omega} \mathfrak{X}$ given by 
\begin{equation*}
f \mapsto \tilde{\mu}_{(\alpha_i), [x_n]_\omega} (f)=\sum_{i,j \; | \; \gamma_i \Lambda = \gamma_j \Lambda \in \Gamma /\Lambda} \alpha_i \alpha_j F (\gamma_i [x_n]_\omega) / (\sum_{i,j \; : \; \gamma_i \Lambda = \gamma_j \Lambda \in \Gamma / \Lambda} \alpha_i \alpha_j),
\end{equation*}
where $F$ is any representative of $f$ in $l^\infty \mathfrak{X}_\omega$.
\end{prop}
\proof 
 Let $\xi= (\xi_n)$ with $\xi_n = \sum_i \alpha_i \delta_{(x_n)} \circ \gamma^{-1}_i$ for every $n$. In virtue of Lemma \ref{lemnew1} we are led to show that the Radon measures $\mu_{\xi/\lim_{n \rightarrow \omega} \|\xi_n\|_2}$ and 
 \begin{equation*}
 \tilde{\mu}_{(\alpha_i), [x_n]_\omega} \circ \pi: f \mapsto \sum_{i,j \; | \; [\gamma_i \Lambda] = [\gamma_j \Lambda] \in \Gamma /\Lambda} \alpha_i \alpha_j f (\gamma_i [x_n]_\omega) / (\sum_{i \in \mathbb{N}} \alpha_i \sum_{j \; | \; \gamma_j \Lambda = \gamma_i \Lambda \in \Gamma / \Lambda} \alpha_j)
 \end{equation*}
  are the same, as states on $l^\infty \mathfrak{X}_\omega$, where $\pi : l^\infty \mathfrak{X}_\omega \rightarrow C(\partial_{\beta, \omega} \mathfrak{X})$ is the quotient map. 
 For, note that for every $n \in \mathbb{N}$ and every $f \in l^\infty \mathfrak{X}$ we have 
 \begin{equation*}
 \begin{split}
  \braket{f \xi_n, \xi_n}&=\sum_{i,j} \alpha_i \alpha_j f (\gamma_i x_n)(\delta_{x_n} \circ \gamma^{-1}_i) (\delta_{x_n} \circ \gamma^{-1}_j)=\sum_i \alpha_i f(\gamma_i x_n)\sum_{j \; : \; \gamma_i x_n = \gamma_j x_n} \alpha_j.
  \end{split}
  \end{equation*}
   Let now $f=(f_n)^\bullet \in l^\infty \mathfrak{X}_\omega$, $\epsilon >0$ and $N \in \mathbb{N}$ be such that $\sum_{i \geq N}\sum_{j \in \mathbb{N}} \alpha_i \alpha_j < \epsilon/(3 \|f\|_\infty)$ and $\sum_{i \in \mathbb{N}} \sum_{j \geq N} \alpha_i \alpha_j < \epsilon /(3\|f\|_\infty)$. There is $A \in \omega$ such that 
   \begin{equation*}
   | \lim_{n\rightarrow \omega}\braket{f_n \xi_n, \xi_n} - \sum_{i} \alpha_i f_n(\gamma_i x_n)\sum_{j : \gamma_i x_n = \gamma_j x_n} \alpha_j| <\epsilon /3
   \end{equation*}
    for every $n \in A$ and so 
   \begin{equation*}
   \begin{split}
   | \lim_{n \rightarrow \omega} \braket{f_n \xi_n, \xi_n} &- \sum_{i \leq N} \alpha_i f_n (\gamma_i x_n)\sum_{j \leq N \; | \; \gamma_i x_n = \gamma_j x_n} \alpha_j| \leq  |  \lim_{n \rightarrow \omega}\braket{f_n \xi_n, \xi_n} - \sum_{i} \alpha_i f_n (\gamma_i x_n)\sum_{j : \gamma_i x_n = \gamma_j x_n} \alpha_j| \\&+  |\sum_{i} \alpha_i f_n (\gamma_i x_n)\sum_{j : \gamma_i x_n = \gamma_j x_n} \alpha_j - \sum_{i\leq N} \alpha_i f_n (\gamma_i x_n)\sum_{j \leq N \; | \; \gamma_i x_n = \gamma_j x_n} \alpha_j| \\ &< | \lim_{n \rightarrow \omega} \braket{f_n \xi_n, \xi_n} - \sum_{i} \alpha_i f_n (\gamma_i x_n)\sum_{j : \gamma_i x_n = \gamma_j x_n} \alpha_j | + 2 \epsilon /3 < \epsilon
   \end{split}
   \end{equation*}
    for every $n \in A$. Similarly, 
    \begin{equation*}
    |\sum_{i \in \mathbb{N}} \alpha_i f(\gamma_i [x_n]_\omega)\sum_{j : \gamma_i \Lambda = \gamma_j \Lambda} \alpha_j - \sum_{i \leq N} \alpha_if(\gamma_i [x_n]_\omega) \sum_{j \leq N : \gamma_i \Lambda = \gamma_j \Lambda} \alpha_j| < \epsilon /3.
    \end{equation*} 
Now, given $(i,j) \in \{0,...,N\} \times \{0,...,N\}$ let $C_{i,j}=\{ n \in \mathbb{N} \; | \; \gamma_i x_n = \gamma_j x_n\}$. Since $\omega$ is an ultrafilter, either $C_{i,j} \in \omega$ or $C_{i,j}^c = \{ n \in \mathbb{N} \; | \; \gamma_i x_n \neq \gamma_j x_n\} \in \omega$. This gives a decomposition of $\{0,...,N\} \times \{0,...,N\}$ in two disjoint sets $\{0,...,N\} \times \{0,...,N\} = D_1 \sqcup D_2$, where $D_1=\{ (i,j) \in \{0,...,N\} \times \{0,...,N\} \; | \; C_{i,j} \in \omega \}$ and $D_2=\{ (i,j) \in \{0,...,N\} \times \{0,...,N\} \; | \; C^c_{i,j} \in \omega\}$. Then taking $D=(\cap_{(i,j) \in D_1} C_{i,j}) \cap (\cap_{(i,j) \in D_2} C^c_{i,j}) \cap A$ we have that $\sum_{i \leq N} \alpha_i f_n (\gamma_i x_n)\sum_{j \leq N : \gamma_i x_n = \gamma_j x_n} \alpha_j= \sum_{(i,j) \in D_1} \alpha_i \alpha_j f_n (\gamma_i x_n)$ for every $n \in D$. Note that $\gamma_i [x_n]_\omega = \gamma_j [x_n]_\omega$ if and only if $(i,j) \in D_1$. Hence $\sum_{i \leq N} \alpha_i f(\gamma_i [x_n]_\omega)\sum_{j \leq N : \gamma_i [x_n]_\omega= \gamma_j [x_n]_\omega} \alpha_j = \sum_{(i,j) \in D_1} \alpha_i \alpha_j f(\gamma_i [x_n]_\omega) $. Now, up to passing to a subset of $D$, say $D' \in \omega$, we have that $|\sum_{(i,j) \in D_1} \alpha_i \alpha_j f(\gamma_i [x_n]_\omega) - \sum_{(i,j) \in D_1} \alpha_i \alpha_j f_n (\gamma_i x_n)| < \epsilon$ for every $n \in D'$. Hence, for $n \in D'$, we have
\begin{equation*}
\begin{split}
| \lim_{n \rightarrow \omega} \braket{f_n \xi_n, \xi_n} &- \sum_{i \in \mathbb{N}} \alpha_i f(\gamma_i [x_n]_\omega)\sum_{j : \gamma_i \Lambda = \gamma_j \Lambda} \alpha_j | \leq | \lim_{n \rightarrow \omega} \braket{f_n \xi_n, \xi_n} - \sum_{i \leq N} \alpha_i f_n (\gamma_i x_n)\sum_{j \leq N \; | \; \gamma_i x_n = \gamma_j x_n} \alpha_j|\\ & + |\sum_{i \leq N} \alpha_i f_n (\gamma_i x_n)\sum_{j \leq N \; | \; \gamma_i x_n = \gamma_j x_n} \alpha_j - \sum_{(i,j) \in D_1} \alpha_i \alpha_j f(\gamma_i [x_n]_\omega)| \\ & + |\sum_{(i,j) \in D_1} \alpha_i \alpha_j f(\gamma_i [x_n]_\omega) - \sum_{i \in \mathbb{N}} \alpha_i f(\gamma_i [x_n]_\omega)\sum_{j : \gamma_i \Lambda = \gamma_j \Lambda} \alpha_j|< 3 \epsilon.
\end{split}
\end{equation*}
It follows that $\lim_{n \rightarrow \omega} \braket{f_n \xi_n, \xi_n} = \sum_{i\in \mathbb{N}} \alpha_i f (\gamma_i [x_n]_\omega) \sum_{j : \gamma_j \Lambda = \gamma_i \Lambda} \alpha_j$. In particular, taking $f=1$ we obtain $(\lim_{n \rightarrow \omega} \|\xi_n\|)^2=\lim_{n \rightarrow \omega} \|\xi_n\|^2  = \sum_{i \in \mathbb{N}} \alpha_i  \sum_{j : \gamma_i \Lambda = \gamma_j \Lambda} \alpha_j$. The result follows. $\Box$

The main reason for specializing to non-trivial free ultrafilters over the natural numbers in our study comes from the fact that the corresponding ultraproducts have a more manageable structure. The following Proposition is the main drawback of this assumption for what concerns dynamics on non-standard boundaries and relies on \cite{HadLi}.

\begin{prop}
\label{prop2}
Let $\Gamma$ be a discrete countable group acting on a countable set $\mathfrak{X}$. Let $\omega \in \partial_\beta \mathbb{N}$ be a non-trivial ultrafilter and $\mu$ be a Calkin measure on $l^\infty \mathfrak{X}_\omega$. Then the image of $l^\infty \mathfrak{X}_\omega$ through the associated GNS-construction equals $L^\infty (\mu)$.
\end{prop}
\proof This is an application of \cite{HadLi} Theorem 4.1, where it is proved that, if we denote by $\pi_\mu$ the GNS-homomorphism, then $\pi_\mu (l^\infty \mathfrak{X}_\omega)$ is a von Neumann algebra and, being weakly dense in $L^\infty (\mu)$, it has to be equal to it. $\Box$\\

The above fact is crucial in order to guarantee the existence of equivariant ucp maps from $l^\infty \Gamma$ to certain inverse limits of $C^*$-algebras. Since $\Gamma$ is exact, such ucp maps give rise to ucp maps from $C^*_\lambda \Gamma$ to the the same inverse limit (extending to $l^\infty \Gamma \rtimes \Gamma$ and then restricting to $C^*_\lambda \Gamma$). In the following we will specialize such construction to the case of certain inverse limits of Koopman representations of $l^\infty \mathfrak{X}_\omega \rtimes \Gamma$ which factor through $C^*(l^\infty \mathfrak{X}_\omega, u_{\sigma, \omega} (\Gamma)) \subset \mathcal{B}_\omega (l^2 \mathfrak{X})$. Before stating the main result, we fix some notation and prove the existence of the desired ucp maps.
Let $A$ be a $C^*$-algebra and consider the preorder on its set of representations given by $\pi \precsim \rho$ if $\pi$ factors through $\rho$, i.e. there is a surjective $*$-homomorphism $\phi_{\pi,\rho} : \rho (A) \rightarrow \pi (A)$ such that $\pi = \phi_{\pi, \rho} \circ \rho$. Let $\Lambda$ be a an upward directed set of representations relative to this preorder, i.e. for every $\pi$, $\rho \in \Lambda$ there is a representation $\theta \in \Lambda$ such that $\pi \precsim \theta$, $\rho \precsim \theta$. If $\pi_\lambda$, $\pi_{\lambda'}$ are elements of $\Lambda$ such that $\pi_\lambda \precsim \pi_{\lambda'}$, we denote by $\pi_{\lambda', \lambda} : \pi_{\lambda'} (A) \rightarrow \pi_\lambda (A)$ the associated surjection. We form the inverse limit of these representations of $A$:
\begin{equation*}
A_\Lambda:= \{ (x_\lambda) \in \prod_\Lambda \pi_\lambda (A) \; | \;  x_\lambda = \pi_{\lambda', \lambda} (x_{\lambda'}) \mbox{ for every } \pi_\lambda \precsim \pi_{\lambda'}\},
\end{equation*}
which is a $C^*$-algebra.

\begin{lem}
\label{lem4}
Let $A$ be a $C^*$-algebra and $A_\Lambda$ an inverse limit as above. The kernel of the $*$-homomorphism $\pi: A \rightarrow A_\Lambda$, $x \mapsto (\pi_\lambda (x))$ equals $\cap_{\lambda \in \Lambda} \ker (\pi_\lambda)$.
\end{lem}
\proof This is clear from the fact that $\pi (x) =0$ if and only if $\pi_\lambda (x) =0$ for every $\lambda \in \Lambda$. $\Box$\\

\begin{lem}
\label{lem3}
Let $A$ be a unital $C^*$-algebra endowed with the action of a discrete countable group $\Gamma$ and let $(\pi_\lambda)_{\lambda \in \Lambda}$ be an upward directed set of representations of $A$ as above, with the additional requirements that $\pi_\lambda (A)$ is a von Neuman algebra admitting an action of $\Gamma$ implemented by weakly-continuous automorphisms for every $\lambda \in \Lambda$, each $\pi_{\lambda, \lambda'}$ is weakly-continuous and $\Gamma$-equivariant for every $\lambda \precsim \lambda' \in \Lambda$. Suppose that for every $\pi_\lambda \in \Lambda$ there is $\Gamma$-equivariant ucp map $\phi_\lambda : l^\infty \Gamma \rightarrow \pi_\lambda (A)$. Then there is a $\Gamma$-equivariant ucp map $\phi : l^\infty \Gamma \rightarrow A_\Lambda$.
\end{lem}
\proof Note that for every finite subset $F \subset \Lambda$, there is a common greater element $\pi_F$ (i.e. $\pi_\lambda \precsim \pi_F$ for every $\pi_\lambda \in F$) and an associated $\Gamma$-equivariant ucp map $\phi_F : l^\infty \Gamma \rightarrow \pi_F (A)$. For each such finite subset $F \subset \Lambda$ consider the $\Gamma$-equivariant ucp map $\tilde{\phi}_F : l^\infty \Gamma \rightarrow \prod_{\lambda \in \Lambda} \pi_\lambda (A)$ given by 
\begin{equation*}
(\tilde{\phi}_F (x) )_\lambda =\begin{cases} \pi_{F, \lambda} (\phi_F (x)) & \pi_\lambda \precsim \pi_F\\
 \phi_\lambda (x) & \pi_\lambda \nprecsim \pi_F, \end{cases}
\end{equation*}
where $\pi_{F, \lambda} : \pi_F (A) \rightarrow \pi_\lambda (A)$ is the $*$-homomorphism given by the condition $\pi_\lambda \precsim \pi_F$. The ucp maps $(\tilde{\phi}_F)$ form a net (indexed by the net of finite subsets of $\Lambda$) and note that the domain of each such map is $l^\infty \Gamma$ and the range is $\prod_{\lambda \in \Lambda} \pi_\lambda (A)$, so they are both von Neumann algebras. Hence, such a net of ucp maps has a limit point in the point-$w^*$-topology (up to passing to a subnet), which we denote by $\phi$. Such limit point takes values in $A_\Lambda$ and is $\Gamma$-equivariant. In order to check the first assertion, represent $\prod_{\lambda \in \Lambda} \pi_\lambda (A)$ faithfully on $\oplus_{\lambda \in \Lambda} H_\lambda$, where for every $\lambda$, $H_\lambda$ is the Hilbert space associated to the representation $\pi_\lambda$, and suppose there are $\lambda_1$, $\lambda_2 \in \Lambda$ such that $\pi_{\lambda_1 , \lambda_2} (z|_{H_{\lambda_2}}) \neq z|_{H_{\lambda_1}}$, where $z=\phi (x)=weak-\lim \phi_i (x)$ for some $x \in A$, where $I$ is the chosen convergent subnet of $\Lambda$ and $i \in I$. There are $\xi, \eta \in H_{\lambda_1}$ and $c >0$ such that $|\braket{(\pi_{\lambda_1 , \lambda_2} (z|_{H_{\lambda_2}}) - z|_{H_{\lambda_1}}) \xi, \eta}| > c$. Note now that $weak-\lim \phi_i (x)|_{H_{\lambda_2}} = z|_{H_{\lambda_2}}$ and $weak-\lim \pi_{\lambda_1, \lambda_2} (\phi_i (x)|_{H_{\lambda_2}})= \pi_{\lambda_1, \lambda_2} (z|_{H_{\lambda_2}})$, since $\pi_{\lambda_1, \lambda_2}$ is weakly-continuous. Hence, given $0 <\epsilon < c/4$ there is $j \in I$ such that $|\braket{(\pi_{\lambda_1, \lambda_2} (\phi_i (x)|_{H_{\lambda_2}}) - \pi_{\lambda_1, \lambda_2} (z|_{H_{\lambda_2}})) \xi, \eta}| < \epsilon$ and $\braket{( \phi_i (x)|_{H_{\lambda_1}} - z|_{H_{\lambda_1}}) \xi, \eta}| < \epsilon$ for every $i \geq j$, entailing $|\braket{(\pi_{\lambda_1 , \lambda_2} (z|_{H_{\lambda_2}}) - z|_{H_{\lambda_1}}) \xi, \eta}| < c/2$, a contradiction. Similarly, $\gamma (z) = \gamma (weak-\lim  \phi_i (x)) = weak-\lim \gamma (\phi_i  x)=weak-\lim (\phi_i (\gamma x))= \phi (\gamma x)$. $\Box$\\

We are now ready for the main result of this work. The idea is to use Lemma \ref{lem3} for a suitable choice of a set $\Lambda$ of Koopman representations in order to obtain a nuclear ucp map $C^*_\lambda \Gamma \rightarrow C^*(l^\infty \mathfrak{X}_\omega, u_{\sigma, \omega} (\Gamma))_\Lambda$; the choice of the representations has to be made in such a way that we are guaranteed that the corresponding ideal is contained in $\mathcal{K}_\omega(l^2 \mathfrak{X})$ for a certain $\omega \in \partial_\beta \mathbb{N}$, a condition which enables the possibility to prove boundary amenability by means of Proposition \ref{prop5} (through Choi-Effros lifting, \cite{ChEf} Theorem 3.10). A technical assumption we use in the following proof is that the group $\Gamma$ satisfies the (AP) property by Haagerup and Kraus (\cite{HaKr}). The assumption that infinite intersections of stabilizer groups are trivial ensures that for elements of the form $a=\sum_{i=1}^k f_{\gamma_i} \gamma_i \in l^\infty \mathfrak{X}_\omega \rtimes \Gamma$, if $\mu$ is a quasi-invariant probability measure of the form appearing in Proposition \ref{prop2.2}, supported on the orbit of a character  $[x_n]_\omega \in {}^\omega \mathfrak{X}_\infty$, and $a$ is in the kernel of the associated Koopman representation, then for every $[y_n]_\omega \in \Gamma [x_n]_\omega$ and every $i=1,...,k$ we have $f_{\gamma_i} ([y_n]_\omega)=0$, since we have $f_{\gamma_i} ([y_n]_\omega) = \lim_{n \rightarrow \omega} f_{\gamma_i, n} (y_n) = \lim_{n \rightarrow \omega} \braket{\sum_{j=1}^k f_{\gamma_j , n} \gamma_j (\gamma_i^{-1} \delta_{y_n}), \delta_{y_n}}_{l^2 \mathfrak{X}}$ and this is realized as a matrix coefficient in the Koopman representation associated to $\mu$. Since ${}^\omega \mathfrak{X}_\infty$ is dense in $\partial_{\beta, \omega} \mathfrak{X}$ (Propositiojn \ref{prop2.1}) this gives $f_{\gamma_i} \in \mathcal{C}_{0,\omega} (\mathfrak{X})$ for every $i$, and so $a \in C^*(C_{0,\omega} (\mathfrak{X}), u_{\sigma,\omega} (\Gamma)) \subset \mathcal{K}_\omega (l^2 \mathfrak{X})$. The (AP)-property simply allows to make similar considerations for general elements in the crossed product.

First we construct a ucp map $\phi : C^*_\lambda \Gamma \rightarrow C^* (l^\infty \mathfrak{X}_\omega, u_{\sigma, \omega} (\Gamma)) \subset \mathcal{B}_\omega (l^2 \mathfrak{X})$ (Proposition \ref{propnew}). Note that the construction of such ucp map only uses the condition that stabilizers of infinite subsets are amenable (rather than trivial). The fact that these are actually trivial, together with the (AP)-property will be used in the proof of Theorem \ref{thm2} in order to guarantee that the kernel of the inverse limit of representations constructed in Proposition \ref{propnew} is contained in $\mathcal{K}_\omega(l^2 \mathfrak{X})$.

\begin{prop}\label{propnew}
Let $\Gamma$ be a discrete countable group acting on a countable set $\mathfrak{X}$ in such a way that the intersection of an infinite number of stabilizer subgroups is always amenable; denote by $\sigma: \Gamma \rightarrow Bij (\mathfrak{X})$ this action. Then there are an inverse limit of representations $\pi$ of $C^*(l^\infty \mathfrak{X}_\omega , u_{\sigma,\omega} (\Gamma))$ associated to Koopman representations on orbits of points in ${}^\omega \mathfrak{X}_\infty$ and a nuclear ucp map $\phi: C^*_\lambda \Gamma \rightarrow \mathcal{B}_\omega(l^2 \mathfrak{X})$ such that for every $\gamma \in \Gamma$ we have $\phi (\lambda(\gamma)) - u_{\sigma, \omega} (\gamma) \in \ker \pi$.
\end{prop}
\proof First, recall from \cite{Ba} Corollary 1.4 that a given an action of $\Gamma$ on a compact Hausdorff space $X$ and a $\Gamma$-quasi-invariant probability measure $\mu$ on $X$, then this measure is Zimmer-amenable if and only the following two conditions are satisfied: (i) the Koopman representation $\pi_\mu$ is tempered, (ii): $\mu$ is pair-amenable, i.e. there is a $\Gamma$-equivariant norm $1$ projection $L^\infty (\mu) \overline{\otimes} L^\infty (\mu) \rightarrow L^\infty (\mu)$. Consider the set of $\Gamma$-quasi-invariant probability measures $\cohull \{\widetilde{NS}_\omega (\mathfrak{X})\}$ as defined in \cite{Ba} Definition 2.1. It follows from Proposition 2.11 of \cite{Ba} that these measures are all pair-amenable. We want to check that these are also tempered. So, first observe that every element $\mu$ in $\cohull \{\widetilde{NS}_\omega (\mathfrak{X})\}$ is a convex combination of uniform approximants of non-standard points $\mu = \sum_{k=1}^l \lambda_k \mu_{(x^{(k)}_n)}$ which are pair-wise singular $\mu_{(x^{(k)}_n)} \perp \mu_{(x^{(k')}_n)}$ for every $k \neq k'$ (the measures $\mu_{(x^{(k)}_n)}$ are constructed from the probability measures associated to the sequence of vectors $(\delta_{x^{(k)}_n})$ and a strictly positive sequence $(\alpha_i) \in l^1 \mathbb{N}$ as in \cite{Ba} Definition 2.1). In virtue of \cite{BaRa} Lemma 2.3 (and since for every finite quasi-invariant measure $\nu$ and for every $\lambda >0$ the Koopman representations on $L^2 (\nu)$ and on $L^2 (\lambda \nu)$ are unitarily equivalent through the unitary $L^2 (\nu) \rightarrow L^2 (\lambda \nu)$, $\xi \mapsto \lambda^{-1/2} \xi$), it is enough to check that every uniform approximant of non-standard points gives rise to a tempered group representation. Let $(x_n)$ be an infinite sequence in $\mathfrak{X}$ and consider the associated point in the spectrum of $l^\infty \mathfrak{X}_\omega$ given by $[x_n]_\omega : f \mapsto \lim_{n\rightarrow \omega} f(x_n)$. Suppose that $\gamma \in \Gamma$ is such that $\gamma ([x_n]_\omega)= [x_n]_\omega$, then there is $A \in \omega$ such that $\gamma (x_n)= x_n$ for every $n \in A$. Since $\omega$ is a free ultrafilter, $|A|=\infty$ and  the stabilizers of these characters are amenable under our assumptions. Note now that, in virtue of Proposition \ref{prop2.2}, for each $k=1,...,l$, the measure $\mu_{(x_n^{(k)})}$ is supported on the orbit $\Gamma [x_n^{(k)}]_\omega = \Gamma / \Gamma_{[x_n^{(k)}]_\omega}$ of $[x_n^{(k)}]_\omega$ (here we view $^{\omega}\mathfrak{X}_\infty$ equivariantly embedded as a dense open discrete subset of $\partial_{\beta, \omega} \mathfrak{X}$ as in Proposition \ref{prop2.1}); since the stabilizer of $[x_n^{(k)}]_\omega$ is amenable, the Koopman representation associated to this measure is thus tempered. It follows that every $\mu \in \cohull \{\tilde{NS}_\omega (\mathfrak{X})\}$ is Zimmer-amenable. It follows from Lemma \ref{lem1} that for every such measure there is a $\Gamma$-equivariant ucp map $l^\infty \Gamma \rightarrow L^\infty (\mu)$. Moreover, for every finite set $\{\mu_1,...,\mu_k\} \subset \cohull \{\widetilde{NS}_\omega (\mathfrak{X})\}$, each of these measures is absolutely continuous with respect to any convex combination of such measures. Hence, considering the associated partial order, every finite set of elements has a common upper bound. We can thus consider $\cohull \{\widetilde{NS}_\omega (\mathfrak{X})\}$ as the index set for a net of representations, namely the Koopman representations of $l^\infty \mathfrak{X}_\omega \rtimes \Gamma$ associated to these measures. Moreover, by Lemma \ref{lemfac} we can consider these as representations of $C^*(l^\infty \mathfrak{X}_\omega , u_{\sigma,\omega} (\Gamma)) \subset \mathcal{B}_\omega(l^2 \mathfrak{X})$. In virtue of Corollary \ref{cor1} we can construct the associated inverse limit $(C^*(l^\infty \mathfrak{X}_\omega , u_{\sigma,\omega} (\Gamma)))_{\cohull \{\widetilde{NS}_\omega (\mathfrak{X})\}}$. In virtue of Proposition \ref{prop2}, Lemma \ref{lem1}, Lemma \ref{lem3}, there is a $\Gamma$-equivariant ucp map $l^\infty \Gamma \rightarrow (l^\infty \mathfrak{X}_\omega)_{\cohull \{\widetilde{NS}_\omega (\mathfrak{X})\}} \subset (C^*(l^\infty \mathfrak{X}_\omega , u_{\sigma,\omega} (\Gamma)))_{\cohull \{\widetilde{NS}_\omega (\mathfrak{X})\}}$, indeed for every $\lambda, \lambda' \in \cohull \{\widetilde{NS}_\omega (\mathfrak{X})\}$ with $\lambda \precsim \lambda'$ the $*$-homomorphism $\pi_{\lambda, \lambda'}|_{L^\infty (\lambda')} : L^\infty (\lambda') \rightarrow L^\infty (\lambda)$ is given by conjugation with the isometry implemented by the Radon-Nikodym derivative as in Proposition \ref{prop1} and so it is weakly-continuous. Since $(C^*(l^\infty \mathfrak{X}_\omega , u_{\sigma,\omega} (\Gamma)))_{\cohull \{\widetilde{NS}_\omega (\mathfrak{X})\}}$ is generated by $(l^\infty \mathfrak{X}_\omega)_{\cohull \{\widetilde{NS}_\omega (\mathfrak{X})\}}$ and the group $\Gamma$ (under the corresponding unitary representation), it follows from \cite{BuEcWi} Lemma 4.8 that we have a ucp map $l^\infty \Gamma \rtimes \Gamma \rightarrow (C^*(l^\infty \mathfrak{X}_\omega , u_{\sigma,\omega}(\Gamma)))_{\cohull \{\widetilde{NS}_\omega (\mathfrak{X})\}}$ which sends group elements to group elements. In particular the homomorphism $\phi: C^*_\lambda \Gamma \rightarrow (C^*(l^\infty \mathfrak{X}_\omega , u_{\sigma,\omega} (\Gamma)))_{\cohull \{\widetilde{NS}_\omega (\mathfrak{X})\}}$ is nuclear and we have the following diagram
\[\begin{tikzcd}
	{C^*_\lambda \Gamma} &&& {(C^*(l^\infty \mathfrak{X}_\omega , u_{\sigma,\omega} (\Gamma)))_{\cohull \{\widetilde{NS}_\omega (\mathfrak{X})\}}} \\
	\\
	&&& {C^*(l^\infty \mathfrak{X}_\omega , u_{\sigma,\omega} (\Gamma))}
	\arrow["\phi"', from=1-1, to=1-4]
	\arrow["\psi", from=1-1, to=3-4]
	\arrow["\pi", from=3-4, to=1-4]
\end{tikzcd}\]
where $\psi$ is given by Choi-Effros extension (\cite{ChEf} Theorem 3.10) and both $\phi$ and $\psi$ are nuclear. Since $\psi$ is an extension of $\phi$ and for every $\gamma \in \Gamma$ we have $\phi(\lambda(\gamma)) = \pi (u_{\sigma, \omega} (\gamma))$, it follows that for every $\gamma \in \Gamma$ we have $\psi (\lambda (\gamma)) - u_{\sigma, \omega} (\gamma) \in \ker \pi$. $\Box$

\begin{thm}
\label{thm2}
Let $\Gamma$ be a discrete countable group with property (AP) (which implies exactness) acting on a countable set $\mathfrak{X}$ in such a way that the intersection of an infinite number of stabilizer subgroups is always trivial. Then the action of $\Gamma$ on $\partial_\beta (\mathfrak{X})$ is topologically amenable.
\end{thm}
\proof 

Let $\pi$ be the representation of $C^*(l^\infty \mathfrak{X}_\omega , u_{\sigma,\omega} (\Gamma))$ given in Proposition \ref{propnew}. In virtue of Proposition \ref{prop5} we want to check that $\ker \pi \subset \mathcal{K}_\omega(l^2 \mathfrak{X})$. \\ Now we will employ property (AP) and triviality of the stabilizers of infinite sets in order to prove the result. In virtue of Lemma \ref{lem4} it is enough to show that $\cap_{\mu \in \cohull \{\widetilde{NS}_\omega (\mathfrak{X})\}} \ker \pi_\mu \subset C^*(\mathcal{C}_{0,\omega} (\mathfrak{X}), u_{\sigma, \omega}(\Gamma))$, where $\mathcal{C}_{0,\omega} (\mathfrak{X}) = l^\infty \mathfrak{X}_\omega \cap \mathcal{K}_\omega(l^2 \mathfrak{X})$. Let $a \in l^\infty \mathfrak{X}_\omega \rtimes \Gamma$. As showed in the proof of \cite{Su} Proposition 3.4 we have $\sigma (a) = \lim_i \sum_{\gamma \in \Gamma} \rho_i (\gamma) f_\gamma \sigma(\gamma)$ (the limit is in norm) for certain finitely supported functions $\rho_i$ on $\Gamma$ which converge pointwise to $1$. Suppose now that $\sigma (a) \in \ker \pi$; we want to show that then $f_\gamma \in \mathcal{C}_{0,\omega} (\mathfrak{X})$ for every $\gamma \in \Gamma$; in order to check this, since $^{\omega}\mathfrak{X}_\infty$ is dense in $\partial_{\beta,\omega} \mathfrak{X}$ (in virtue of Proposition \ref{prop2.1}), it will be enough to check that $f_\gamma ([x_n]_\omega) =0$ for every $\gamma \in \Gamma$ and every $[x_n]_\omega \in {}^{\omega}\mathfrak{X}_\infty$. Suppose then that there are an infinite sequence $(y_n) \subset \mathfrak{X}$ and a $\gamma \in \Gamma$ such that $\lim_{n \rightarrow \omega}f_\gamma (y_n)= \alpha \neq 0$ (in virtue of Proposition \ref{prop2.1} this is equivalent to the fact that $f_\gamma \notin \mathcal{C}_{0,\omega} (\mathfrak{X})$). Let $i$ be large enough so that $\rho_i (\gamma) > 1/2$. Consider a uniform approximant of non-standard points $\mu$ associated to $(x_n)=(\gamma^{-1} y_n)$ and the corresponding sequence of $l^2$-functions $(\xi_n)$, as in the paragraph preceding Proposition \ref{prop1.5}. As in the last part of the proof of \cite{BaRa} Theorem 2.7 (or equivalently Theorem \ref{thmap2} in the appendix) we see that the linear functional on $\mathbb{C}[\Gamma]$ given by $\eta \mapsto \lim_{n \rightarrow \omega} \braket{\delta_{\eta x_n}, \delta_{\gamma x_n}}_{l^2 \mathfrak{X}}$ is implemented by vectors in $L^2 (\mu)$, namely, for every $\eta \in \Gamma$ we have $\lim_{n \rightarrow \omega} \braket{\delta_{\eta x_n}, \delta_{\gamma x_n}}_{l^2 \mathfrak{X}} = \braket{\pi_\mu (\eta) (\delta_{x_n}/\xi_n), (\delta_{\gamma x_n}/\xi_n)}_{L^2 (\mu)}$. For every $j \geq i$ we have
\begin{equation*}
\begin{split}\braket{\pi_\mu^\rtimes (\sum_{\eta \in \Gamma} \rho_j (\eta) f_\eta \eta)\delta_{ x_n}/\xi_n, \delta_{\gamma x_n}/\xi_n}_{L^2(\mu)} &=\lim_{n \rightarrow \omega} \braket{\sum_{\eta \in \Gamma} \rho_j (\eta) f_\eta \delta_{\eta x_n}, \delta_{\gamma x_n}}_{l^2 \mathfrak{X}}\\
&= \lim_{n \rightarrow \omega} \rho_j (\gamma) f_\gamma (\gamma x_n)=\lim_{n \rightarrow \omega} \rho_j (\gamma) f_\gamma (y_n),
\end{split}
\end{equation*}
whose absolute value is greater than $|\alpha|/2$, where we used the hypothesis that stabilizers of infinite subsets are trivial in order to deduce that $\lim_{n \rightarrow \omega} \braket{ \delta_{\eta x_n}, \delta_{\gamma x_n}}_{l^2 \mathfrak{X}}= \delta_{\gamma, \eta}$ for every $\eta \in \Gamma$. In particular $\pi_\mu^\rtimes (a) \neq 0$ and $a \notin \ker \pi_\mu^\rtimes$. The result follows. $\Box$\\

\begin{oss} The assumption that $\Gamma$ has property (AP) is in general not essential. As an example, the action of $\SL(3,\mathbb{Z})$ on the coset space $\SL(3,\mathbb{Z}) / \SL(2,\mathbb{Z})$ (where we view $\SL(2,\mathbb{Z}) < \SL(3,\mathbb{Z})$ as the subgroup of matrices $(a_{i,j})_{i,j=1}^3$ with $a_{1,2}=a_{1,3}=a_{21}=a_{3,1}=0$) is certainly not amenable, but is amenable when we pass to the Stone-{\v C}ech boundary $\partial_\beta (\SL(3,\mathbb{Z}) / \SL(2,\mathbb{Z}))$ (\cite{BaRa2}). The fact that $\SL(3,\mathbb{Z})$ does not have property (AP) is proved in \cite{LaSa}. \end{oss}

\begin{oss} In general the conclusion of Theorem \ref{thm2} does not hold if we relax the hypothesis on the infinite intersection of stabilizer subgroups by assuming it to be amenable, rather than trivial (in particular the kernel of the homomorphism $\pi$ constructed in Proposition \ref{propnew} is not always contained in $\mathcal{K}_\omega (l^2 \mathfrak{X})$). Indeed, let $\Gamma = \SL(2, \mathbb{Z}[1/p])$ and consider the left-right action of $\Gamma \times \Gamma$ on $\Gamma$. In virtue of \cite{BoPi} Corollary 4.3 $\SL(2,\mathbb{Q}_p)$ is weakly amenable, hence has property (AP); thus, by \cite{HaKr} Theorem 2.4 and \cite{HaKr} Corollary 1.17 $\Gamma \times \Gamma$ has the (AP)-property as well. Note also that the infinite intersection of stabilizer subgroups is always abelian. In order to see this, first observe that the centralizer of an element $\gamma$ in $\SL(2, \mathbb{Z}[1/p])$ is abelian, actually either conjugated to diagonal matrices (if $\gamma$ has two distinct eigenvalues) or is conjugated to strictly upper-triangular matrices (if $\gamma$ has one eigenvalue). Hence suppose $\gamma, \eta, x, y \in \Gamma$ are such that $x \neq y$ and $\gamma x \eta^{-1} = x$, $\gamma y \eta^{-1} = y$. Then $\gamma xy^{-1} \gamma^{-1} = xy^{-1}$ and $\eta x^{-1} y \eta^{-1} = x^{-1} y$ and so both $\gamma$ and $\eta$ belong to an abelian group. Hence infinite intersections of stabilizer subgroups for the action of $\Gamma \times \Gamma$ on $\Gamma$ are abelian, hence amenable. But the action of $\Gamma \times \Gamma$ on $\partial_\beta \Gamma$ is not amenable, since $\Gamma$ is a lattice in $\SL(2,\mathbb{R}) \times \SL(2,\mathbb{Q}_p)$ and as such is not bi-exact (by \cite{De} Theorem E).
\end{oss}


The following is a folklore application of Theorem \ref{thm2}. 

\begin{cor}
Let $\Gamma$ be a discrete countable group with property (AP) acting on a countable set $\mathfrak{X}$. If the induced action on $\partial_\beta \mathfrak{X}$ is free, then it is amenable.
\end{cor}
\proof This is a consequence of the fact that for every $\omega \in \partial_\beta \mathbb{N}$, the non-standard boundary $\partial_{\beta, \omega} \mathfrak{X}$ is a $\Gamma$-equivariant extension of $\partial_\beta \mathfrak{X}$ which contains a copy of ${}^\omega \mathfrak{X}_\infty$. $\Box$

\section{Boundary amenable, non-amenable dynamical systems}
This section contains the claimed examples of non-amenable, boundary amenable dynamical systems as an application of Theorem \ref{thm2}.

\begin{thm}
\label{thm1}
Let $\Gamma$ be a countable hyperbolic group without torsion and $\Lambda$ be an infinite quasi-isometrically embedded subgroup with infinite index. Then the action of $\Gamma$ on $\partial_\beta (\Gamma / \Lambda)$ is topologically amenable. In particular, if $\Lambda$ is a finitely generated subgroup of a free group $\mathbb{F}_r$, $r \geq 2$, and $\Lambda$ has infinite index in $\mathbb{F}_r$, then the action of $\mathbb{F}_r$ on $\partial_\beta (\mathbb{F}_r / \Lambda)$ is topologically amenable.
\end{thm}
\proof Note that by \cite{Ozh} hyperbolic groups are weakly amenable, hence they have the (AP)-property. In virtue of \cite{BrHa} Corollary III.$\Gamma$.3.6 $\Lambda$ is quasi-convex in $\Gamma$ if and only if it is quasi-isometrically embvedded. The result is then an application of Theorem \ref{thm2} and \cite{Gi} Theorem 2.2, where it is proved that quasi-convex subgroups of hyperbolic groups have finite height and, since $\Gamma$ is torsion free, this implies that the intersection of an infinite number of conjugates of $\Lambda$ is always trivial. The fact that every finitely generated subgroup of a free group is quasi-convex is well known and is a direct consequence of the fact that the Cailey graph of a free group is a tree. $\Box$


\begin{thm}
Let $\Gamma$ be a discrete countable group and $\Lambda$ be a discrete countable infinite group. Then the action of $\Upsilon$ on $\partial_\beta (\Gamma * \Lambda / \Gamma)$ is amenable for every subgroup $\Upsilon$ of $\Gamma * \Lambda$ with the (AP). Moreover, if $\Gamma$ is non-amenable and $\Gamma * \Lambda$ is $C^*$-simple and has property (AP), the only non-trivial ideal in $\pi_{\Gamma * \Lambda/\Gamma}(C^*(\Gamma * \Lambda))$ is the ideal of compact operators.  In particular $\pi_{\lambda_{\mathbb{F}_r * \mathbb{F}_s / \mathbb{F}_r} }(C^* \mathbb{F}_{r+s})$ has a unique ideal for every $r,s \geq 2$.
\end{thm}
\proof Let $\eta \Gamma$ and $\eta' \Gamma$ be different cosets in $\Gamma * \Lambda / \Gamma$. Then an element $s \in \Gamma * \Lambda$ fixes both these sets if and only if $s \in \eta' \Gamma (\eta')^{-1} \cap \eta \Gamma \eta^{-1} = \{e\}$. Hence the first statement follows from Theorem \ref{thm2}. For the second statement, note that $\Gamma$ is a-normal in $\Gamma * \Lambda$ and so has spectral gap (see \cite{BeKa} Section 4), hence by \cite{BeKa} Theorem B $\pi_{\Gamma * \Lambda/\Gamma}(C^*(\Gamma * \Lambda))$ contains the compact operators as a minimal ideal. Since the action of $\Gamma * \Lambda$ on $\partial_\beta (\Gamma * \Lambda/\Gamma)$ is topologically amenable, it follows that $\pi_{\Gamma * \Lambda/\Gamma}(C^*(\Gamma * \Lambda)) / \mathbb{K}(l^2(\Gamma * \Lambda/\Gamma)) \simeq C^*_\lambda (\Gamma * \Lambda)$ and so $\mathbb{K}(l^2(\Gamma * \Lambda/\Gamma))$ is the only non-trivial ideal in $\pi_{\Gamma * \Lambda/\Gamma}(C^*(\Gamma * \Lambda))$. $\Box$

\begin{cor}
Let $\Lambda$ and $\Gamma$ be infinite countable groups which are weakly amenable with Cowling-Haagerup constant $\Lambda=1$ and suppose $\Gamma$ is not amenable. Then the action of $\Gamma * \Lambda$ on $\partial_\beta (\Gamma * \Lambda / \Gamma)$ is topologically amenable and the associated $C^*$-algebra $\pi_{\Gamma * \Lambda/\Gamma}(C^*(\Gamma * \Lambda))$ has a unique non-trivial ideal.
\end{cor}
\proof This follows from the fact that the class of weakly amenable groups with Cowling-Haagerup constant $\Lambda=1$ is stable under free products (\cite{RiXu}) and that weakly amenable groups have property (AP). $\Box$

Other examples of actions of discrete countable groups on countable sets for which stabilizers of infinite sequences are always trivial have been considered in \cite{AbGl}, where it is proved that if $\Gamma$ is a countable group acting on a $k$-regular tree $T$, with $k \geq 3$, in such a way that the infinite intersection of stabilizer subgroups is always trivial, then for almost (in a suitable sense) all elements $\gamma \in \Aut(T)$, the group generated by $\Gamma$ and $\gamma$ shares the same property (\cite{AbGl} Theorem). Hence we have the following

\begin{thm}
Consider the topology on the group $\Aut (T)$ of automorphisms of a $k$-regular tree $T$ ($k \geq 3$) given by pointwise convergence and the associated Haar measure, normalized so that vertex stabilizers have measure $1$. Let $a_1,..., a_n$ be independent Haar-random elements of $\Aut (T)$ and let $\Gamma = \braket{a_1,...,a_n}$. Then almost surely, the action of $\Gamma$ on $\partial_\beta T$ is topologically amenable.
\end{thm}
\proof This is a direct application of the Corollary in \cite{AbGl}. Indeed, if $a_1,..., a_n$ are independent Haar-random elements of $\Aut(T)$, then almost surely the group generated by them is a free group of rank $n$ which acts on $T$ in such a way that the infinite intersection of stabilizer subgroups is trivial. The result follows then from Theorem \ref{thm2}. $\Box$

\section{Appendix: Dynamics on non-standard boundaries}

This appendix contains results (with proof sketches) mainly taken from \cite{Ba} which are used throughout the previous sections.\\ 
In \cite{Ba} C. Anantharaman Delaroche gave a characterization of Zimmer-amenability for actions of second countable groups on standard spaces. The same characterization can be obtained in the case of actions of discrete countable groups on non-separable spaces. In fact we have the following
\begin{thm}
Let $\Gamma$ be a discrete countable group acting on a compact Hausdorff space $X$. Let $\mu$ be a $\Gamma$-quasi-invariant probability measure on $X$. Then $\mu$ is Zimmer-amenable if and only if the following two conditions are satisfied:
\begin{itemize}
\item[(i)] the associated Koopman representation $\pi_\mu$ is tempered,
\item[(ii)] there is a $\Gamma$-equivariant norm-$1$ projection $L^\infty (\mu) \bar{\otimes} L^\infty (\mu) \rightarrow L^\infty (\mu)$, i.e. $\mu$ is pair-amenable.
\end{itemize}
\end{thm}
\proof[Sketch of proof] First of all the fact that Zimmer-ameanbility of $\mu$ implies properties $(i)$ and $(ii)$ is really the same argument appearing in \cite{AnDe2} Proposition 4.3.2. So we only need to check that if $(i)$ and $(ii)$ are satisfied then $\mu$ is Zimmer-amenable. In order to see this, first of all we employ arguments similar to the ones appearing in \cite{BeCr} in order to see that the existence of a $\Gamma$-equivariant norm-$1$ projection $L^\infty (\mu) \bar{\otimes} L^\infty (\mu) \rightarrow L^\infty (\mu)$ is equivalent to the existence of approximately invariant vectors in $L^2 ((\mu), L^\infty (\mu))$ (in a suitable sense, see \cite{Ba} Proposition 2.1); then, since $\pi_\mu$ is tempered, we can employ Voiculescu's Theorem (\cite{Da} Theorem II.5.3) to produce an isometry $V: l^2 \Gamma \rightarrow L^2 (\mu)$. Then, if $\xi$ is $\epsilon/3$-approximately invariant in $L^2 (\mu, L^\infty (\mu))$ (in the sense given at the beginning of the proof of Theorem 1.3 in \cite{Ba}), it follows that $(V \otimes 1) \xi$ is $\epsilon$-approximately invariant in $l^2 (\Gamma, L^\infty (\mu))$ and so at the limit we obtain an equivariant norm-$1$ projection $l^\infty (\Gamma) \bar{\otimes} L^\infty (\mu) \rightarrow L^\infty (\mu)$. $\Box$


The next result explains the relationship between regularity properties of Koopman representations associated to Calkin measures and regularity properties of the Calkin representation associated to an action of a discrete countable group on a countable set.

\begin{defn}
Let $\Gamma$ be a discrete countable group and let $C^*_\tau \Gamma$ be a $C^*$-completion of the group algebra $\mathbb{C}[\Gamma]$. A map from $C^* \Gamma$ to a set $S$ is said to be \textit{$\tau$-continuous} if it factors through the canonical surjection $C^* \Gamma \rightarrow C^*_\tau \Gamma$.
\end{defn}

\begin{thm}
\label{thmap2}
Let $\Gamma$ be a discrete countable group acting on a countable set $\mathfrak{X}$ by means of a group homomorphism $\sigma : \Gamma \rightarrow Bij (\mathfrak{X})$ and let $C^*_\tau \Gamma$ be a $C^*$-completion of the group algebra $\mathbb{C}[\Gamma]$. Then the associated Calkin representation $Q_\mathfrak{X} \circ u_\sigma$ is $\tau$-continuous if and only if for every $\omega \in \partial_\beta \mathbb{N}$ and every Calkin measure $\mu$ on $\partial_{\beta, \omega} (\mathfrak{X})$, the associated Koopman representation $\pi_\mu$ is $\tau$-continuous.
\end{thm}
\proof[Sketch of proof] Suppose that the Koopman representations associated to Calkin measures are all $\tau$-continuous. Note that the requirement that the Calkin representation $Q_{\mathfrak{X}} \circ u_\sigma$ is $\tau$-continuous is equivalent to the requirement that every state $\phi$ on $\mathcal{Q}(l^2 \mathfrak{X})$ is $\tau$-continuous. Now, states on $\mathcal{Q}(l^2 \mathfrak{X})$ are all given by pairs $(\omega, (\xi_n))$, where $\omega \in \partial_\beta \mathbb{N}$ is a free ultrafilter on the natural numbers and $(\xi_n)$ is a sequence of functions in $l^2 \mathfrak{X}$ of norm one such that $\xi_n \rightarrow 0$ weakly with respect to $\omega$; explicitely every such $\phi$ has the form $\phi (a)= \lim_{n \rightarrow \omega} \braket{\bar{a} \xi_n,\xi_n}_{l^2 \mathfrak{X}}$, for $a \in \mathcal{Q}(l^2 \mathfrak{X})$ and $\bar{a} \in \mathbb{B}(l^2 \mathfrak{X})$ such that $Q_{\mathfrak{X}} (\bar{a})=a$. So the idea is to realize that given such a state, when restricted to $u_\sigma (C^* \Gamma)$, it can be realized as a vector state for the Koopman representation associated to an appropriate Calkin measure on the non-standard boundary $\partial_{\beta, \omega} \mathfrak{X}$ (the ultrafilter appearing here is the same defining the state on the Calkin algebra). So, first of all, we construct the desired Calkin measure: for, let $(\alpha_i)_{i \in \mathbb{N}} \in l^1 \mathbb{N}$ be such that $\alpha_i >0$ for every $i \in \mathbb{N}$ and let $\Gamma = \{\gamma_i\}_{i \in \mathbb{N}}$ be an enumeration of $\Gamma$ with $\gamma_0 = \id$. Define, for every $n \in \mathbb{N}$, the vector $\hat{\xi}_n := \sum_i \alpha_i | u_\sigma (\gamma_i) (\xi_n)|$ and normalize it to obtain an $l^2$-function of norm one, $\tilde{\xi}_n := \hat{\xi}_n / \| \hat{\xi}_n\|_2$. The probability measure on $\sigma(l^\infty \mathfrak{X}_\omega)$ given by $f \mapsto \lim_{n \rightarrow \omega} \braket{f \tilde{\xi}_n , \tilde{\xi}_n}_{l^2 \mathfrak{X}}$ is then supported on $\partial_{\beta, \omega} \mathfrak{X}$ (since $\xi_n \rightarrow 0$ weakly with respect to $\omega$) and is the desired Calkin measure $\mu_\phi$. Now we want to check that the original state $\phi$ is implemented by a vector in the associated Koopman representation. In order to do so, we use the characterization of the Koopman representation among the group homomorphisms $\Gamma \rightarrow U L^2 (\mu_\phi)$ established in \cite{BaRa} Lemma 2.5, i.e. $\pi_{\mu_\phi}$ is the only such group homomorphism which gives a covariant pair $(\pi_{\mu_\phi}, M)$, where $M: C(\partial{\beta,\omega} \mathfrak{X}) \rightarrow \mathbb{B}(L^2 (\mu_\phi))$ is the representation as multiplication operators, and which sends positive functions in $L^2 (\mu_\phi)$ to positive functions. Given such characterization we exhibit a group homomorphism satisfying the above two properties in the following way: first consider the linear subspace of $\prod_{\mathbb{N}} l^\infty \mathfrak{X}$ given by $W:=\{ (f_n) \; | \; \forall \gamma \in \Gamma \; \exists \; M_\gamma >0 \;, \; |u_\sigma (f_n) | \leq M_\gamma \tilde{\xi}_n \; \forall n\}$ and let $\tilde{W}:= \{ (h_n) \; | \; h_n = f_n / \tilde{\xi}_n, \; (f_n) \in W\}$. This is a dense subspace of $L^2 (\mu_\phi)$, which is invariant under $M(C(\partial_{\beta,\omega} \mathfrak{X}))$. For every $\gamma \in \Gamma$ the map $\pi_\gamma : \tilde{W} \rightarrow \tilde{W}$, $f_n /\tilde{\xi}_n \mapsto u_\sigma (\gamma)(f_n) /\tilde{\xi}_n$ satisfies the assumptions of the Lemma, when restricted to $\tilde{W}$; from this it is possible to see that also the extension to $L^2 (\mu_\phi)$ satisfies such properties. Then the matrix coefficient associated to $\xi_n /\tilde{\xi}_n$ gives the state $\phi$, indeed, for every $\gamma \in \Gamma$ we have $\phi(u_\sigma (\gamma)) = \lim_{n \rightarrow \omega} \braket{u_\sigma (\gamma) \xi_n , \xi_n}_{l^2 \mathfrak{X}} = \braket{\pi_{\mu_\phi} (\gamma) (\xi_n /\tilde{\xi}_n), (\xi_n /\tilde{\xi}_n)}_{L^2 (\mu_\phi)}$.\\
For the converse, we use the fact that the vector states associated to elements in $\tilde{W}$ all give rise to states on $\mathcal{Q}(l^2 \mathfrak{X})$ and the fact that $\tilde{W}$ is dense in $L^2 (\mu_\phi)$. $\Box$

Clearly, if a discrete countable group $\Gamma$ acts on a countable set $\mathfrak{X}$ in such a way that the induced action on $\partial_\beta \mathfrak{X}$ is amenable, then for every $\omega \in \partial_\beta \mathbb{N}$ the action on the associated non-standard boundary $\partial_{\beta, \omega} \mathfrak{X}$ is amenable as well, since these are all equivariant extensions of $\partial_\beta \mathfrak{X}$. Notably, each of these extensions already capture amenability of the action of $\Gamma$ on $\partial_\beta \mathfrak{X}$, in fact the following holds:
\begin{prop}[\cite{Ba} Proposition 2.3]
\label{prop5}
Let $\mathfrak{X}$ be a countable set endowed with an action $\sigma$ of a countable exact group $\Gamma$ and denote by $u_\sigma$ the associated unitary representation on $l^2 \mathfrak{X}$. Let $\omega \in \partial_\beta \mathbb{N}$ and denote by $u_{\sigma, \omega}$ the unitary group representation of $\Gamma$ in $\mathcal{B}_\omega (l^2 \mathfrak{X})$ given by $u_{\sigma, \omega} (\gamma)= (u_\sigma (\gamma))$. The following are equivalent:
\begin{itemize}
\item[(i)] The action of $\Gamma$ on $\partial_{\beta, \omega} \mathfrak{X}$ is topologically amenable,
\item[(ii)] there is a nuclear ucp map $\phi: C^*_\lambda \Gamma \rightarrow \mathcal{B}_\omega(l^2 \mathfrak{X})$ such that $\phi (\lambda (\gamma)) - u_{\sigma, \omega} (\gamma) \in \mathcal{K}_\omega (l^2 \mathfrak{X})$ for every $\gamma \in \Gamma$,
\item[(iii)] the action of $\Gamma$ on $\partial_\beta \mathfrak{X}$ is topologically amenable.
\end{itemize}
\end{prop}
\proof[Sketch of proof] Clearly $(iii) \Rightarrow (i)$ since $\partial_{\beta,\omega} \mathfrak{X}$ is a $\Gamma$-equivariant extension of $\partial_\beta \mathfrak{X}$. To see that $(i) \Rightarrow (ii)$ just note that, by amenability of the action of $\Gamma$ on $\partial_{\beta, \omega} \mathfrak{X}$, if we denote by $\iota_\omega : C(\partial_{\beta, \omega} \mathfrak{X}) \rightarrow \mathcal{B}_\omega(l^2 \mathfrak{X})$ the canonical embedding, the covariant pair $( \iota_\omega, u_{\sigma, \omega})$ gives rise to a representation of the reduced crossed product $C(\partial_{\beta, \omega} \mathfrak{X})\rtimes_\lambda \Gamma$, since this is nuclear, we obtain, by restriction to $C^*_\lambda \Gamma$, the desired nuclear ucp map. For the implication $(ii) \Rightarrow (iii)$, first take a nuclear ucp lift of $\phi$, say $ \psi: C^*_\lambda \Gamma \rightarrow \prod_\mathbb{N} \mathbb{B}(l^2 \mathfrak{X})$. Let $F_k \subset F_{k+1} \subset \Gamma$ be an exhaustion of $\Gamma$ by finite symmetric sets. Since $\psi$ is a ucp lift of $\phi$, there are sequences $n_k \rightarrow \infty$, $A_k \in \omega$ nested sets such that for every $\gamma \in F_k$ we have
$\| P_{n_k} (\psi_i (\lambda (\gamma)) - u_{\sigma} (\gamma))P_{n_k} - (\psi_i (\lambda(\gamma)) - u_\sigma (\gamma))\| < 1/k^2$ for every $i \in A_k$ ($\psi_i$ is the projection
on the $i$-th factor of $\psi$). Choosing for every $k$ any $i \in A_k$ one finds isometries $V_k : l^2 \mathfrak{X} \rightarrow l^2 \Gamma$ by means of Voiculescu's Theorem (\cite{Da} Theorem II.5.3). Proceeding as in the proof of \cite{Oz4} Proposition .4.1 we construct maps $\mu_k : \mathfrak{X} \rightarrow \mathcal{P}(\Gamma)$ which are approximately equivariant outside certain finite subsets of $\mathfrak{X}$, taking an appropriate average of these maps we obtain a map $\mu : \mathfrak{X} \rightarrow \mathcal{P}(\Gamma)$ which satisfies $\lim_{x \rightarrow \infty} \| \gamma \mu (x) - \mu (\gamma x)\|_1 =0$ for every $\gamma \in \Gamma$ and the conclusion follows from \cite{Oz4} Proposition 4.1. $\Box$

It follows that the action of $\Gamma$ on $\partial_\beta \mathfrak{X}$ is topologically amenable if and only if the same is true for the action of $\Gamma$ on a non-standard boundary $\partial_{\beta, \omega} \mathfrak{X}$ for some $\omega \in \partial_\beta \mathbb{N}$ (and hence for every $\omega \in \partial_\beta \mathbb{N}$); this amounts to the fact that every quasi-invariant probability measure on $\partial_{\beta, \omega} \mathfrak{X}$ is Zimmer-amenable. By the above characterization we than want to understand for which measures pair-amenability is automatic. The procedure used in the proof of Theorem \ref{thmap2} to construct measures from sequences of vectors in $l^2 \mathfrak{X}$ give the desired property in the case the starting vector is of the form $\xi_n = \delta_{x_n}$ for a certain sequence $(x_n) \subset \mathfrak{X}$. Measures of this form are referred to as "uniform approximants of non-standard points" in \cite{Ba}. In fact, as shown in \cite{Ba} Proposition 2.11, the following holds:
\begin{prop}
Let $\mu$ be a quasi-invariant Radon probability measure on $\sigma (l^\infty \mathfrak{X}_\omega)$ which is a convex combination of uniform approximants of non-standard points associated to sequences $(\alpha_i^{(m)})$ satisfying $\alpha_i^{(m)} >0$ for every $i,m$. The pair $(\sigma (l^\infty \mathfrak{X}_\omega) \times \sigma (l^\infty \mathfrak{X}_\omega), \mu \otimes \mu), (\sigma (l^\infty \mathfrak{X}_\omega), \mu))$ is amenable, i.e. there is an equivariant norm-$1$ projection $L^\infty (\mu) \bar{\otimes} L^\infty (\mu) \rightarrow L^\infty (\mu)$.
\end{prop}
The main ingredient in the proof of the above result is that convex combinations of uniform approximants of non-standard points are measures induced by sequences $(\xi_n) \subset l^2 \mathfrak{X}$ which do not give rise to trivial classes in $l^\infty \mathfrak{X}_\omega$. This gives the possibility to construct an isometry $L^2( \mu) \rightarrow L^2 (\mu )\otimes L^2( \mu)$ which implements the desired equivariant projection $L^\infty (\mu ) \otimes L^\infty (\mu) \rightarrow L^\infty (\mu)$.

\section{Acknowledgments}
The author acknowledges the support of the grant SOE-Young Researchers 2024 "Generalized Akemann-Ostrand Property:
Analytical, Dynamical And Rigidity Properties (GAOPADRP)",
CUP: E83C24002550001 and the MIUR Excellence Department Project
2023--2027 MatMod@Tov awarded to the Department of Mathematics,
University of Rome Tor Vergata. He also acknowledges the support of INdAM-GNAMPA and of IM PAN, where this work was partially developed. The author thanks the anonymous referee for the useful suggestions which led to an improvement of the exposition and of the range of applications of the results proved in this work.


\baselineskip0pt
\bigskip
  \footnotesize


\end{document}